\documentclass[a4paper,11pt,reqno]{amsart}
\usepackage{graphicx}
\usepackage{amsmath}
\usepackage{amssymb}
\usepackage{oldgerm}
\usepackage[all]{xy}
\newtheorem{te}{Theorem}[section]
\newtheorem{co}[te]{Corollary}
\newtheorem{de}[te]{\sc Definition}
\newtheorem{ex}[te]{Example}
\newtheorem{prop}[te]{Proposition}
\newtheorem{lm}[te]{Lemma}
\begin{document}
\baselineskip=15pt
\title{Graded Brauer tree algebras}
\author{Dusko Bogdanic}
\address{Dusko Bogdanic \newline Mathematical Institute \\ University of Oxford \\  \newline 24-29 St.\
Giles \\ Oxford OX1 3LB \\ United Kingdom}

\email{bogdanic@maths.ox.ac.uk}
\date{}

\begin{abstract}
In this paper we construct non-negative gradings on a basic Brauer
tree algebra $A_{\Gamma}$ corresponding to an arbitrary Brauer
tree $\Gamma$ of type $(m,e)$. We do this by transferring gradings
via derived equivalence from a basic Brauer tree algebra $A_S$,
whose tree is a star with the exceptional vertex in the middle, to
$A_{\Gamma}$. The grading on $A_S$ comes from the tight grading
given by the radical filtration. To transfer gradings via derived
equivalence we use tilting complexes constructed by taking Green's
walk around $\Gamma$ (cf.\ [\ref{Zak}]). By computing endomorphism
rings of these tilting complexes we get graded algebras.

We also compute ${\rm Out}^K(A_{\Gamma})$, the group of outer
automorphisms that fix isomorphism classes of simple
$A_{\Gamma}$-modules, where $\Gamma$ is an arbitrary Brauer tree,
and we prove that there is unique grading on $A_{\Gamma}$ up to
graded Morita equivalence and rescaling.
\end{abstract}
\maketitle \vspace{-10mm}
\section{Introduction}
\noindent In this paper we transfer gradings between Brauer tree
algebras via derived equivalences. Our work has been motivated by
Theorem 6.4 in an unpublished paper [\ref{Rou}] by Rouquier. This
theorem says that gradings are compatible with derived
equivalences. We show how this idea is used on the class of Brauer
tree algebras. For an arbitrary Brauer tree $\Gamma$ of type
$(m,e)$, ie.\ for a Brauer tree with $e$ edges and multiplicity of
the exceptional vertex $m$, we transfer tight grading from the
basic Brauer tree algebra $A_S$ corresponding to the Brauer star
$S$ of the same type $(m,e)$, to the Brauer tree algebra
$A_{\Gamma}$. In Section 4 we prove that the resulting grading on
$A_{\Gamma}$ is non-negative and we investigate its properties. In
particular, this construction associates to each Brauer tree
algebra $A$, which is a symmetric algebra, a quasi-hereditary
algebra $A_0$, the subalgebra of $A$ consisting of the elements of
$A$ which have degree 0. We prove in Section 5 that the knowledge
of the subalgebra $A_0$ and of the cyclic ordering of its
components is sufficient to recover the whole algebra $A$. In
Sections 6 and 7 we give explicit formulae for the graded Cartan
matrix  and graded Cartan determinant of $A_{\Gamma}$, and we
prove that the graded Cartan determinant only depends on the type
of the Brauer tree. Sections 9 and 10 deal respectively with the
problem of shifting summands of a tilting complex and with the
change of the exceptional vertex when the multiplicity of the
exceptional vertex is 1. In the last section we compute ${\rm
Out}^K(A_{\Gamma})$, the group of outer automorphisms that fix
isomorphism classes of simple $A_{\Gamma}$-modules, and we prove
that it only depends on the multiplicity of the exceptional
vertex. We also classify all gradings on an arbitrary Brauer tree
algebra, and we prove that there is unique grading up to graded
Morita equivalence and rescaling.

\section{Notation} Throughout this text $k$ will be an
algebraically closed field of positive characteristic. All
algebras will be finite dimensional algebras over $k$ and all
modules will be left modules. The category of finite dimensional
$A$--modules will be denoted by $A$--${\rm mod}$ and the full
subcategory of finite dimensional projective $A$--modules will be
denoted by $P_A$. The derived category of bounded complexes over
$A$--${\rm mod}$ will be denoted by $D^b(A)$ and the homotopy
category of bounded complexes over $P_A$ will be denoted by
$K^b(P_A)$.

Let $A$ be a $k$--algebra. We say that $A$ is a graded algebra if
$A$ is  the direct sum of subspaces $A=\bigoplus_{i\in\mathbb{Z}}
A_i$, such that $A_iA_j\subset A_{i+j}$, $i,j\in \mathbb{Z}$. If
$A_i=0$  for $i<  0$, we say that $A$ is non-negatively graded. An
$A$-module $M$ is graded if it is the direct sum of its subspaces
$M=\bigoplus_{i\in\mathbb{Z}} M_i$,  such that  $A_iM_j\subset
M_{i+j}$, for all $i,j\in \mathbb{Z}$. If $M$ is a graded
$A$--module, then $N=M\langle i\rangle$ denotes the graded module
given by $N_j=M_{i+j}$, $j\in \mathbb{Z}$. An $A$-module
homomorphism $f$ between two graded modules $M$ and $N$ is a
homomorphism of graded modules if $f(M_i)\subseteq N_i$, for all
$i\in \mathbb{Z}$. For a graded algebra $A$, we denote by
$A$--${\rm modgr}$ the category of graded $A$--modules. This
category sits inside the category of graded vector spaces. Its
objects are graded $A$--modules and morphisms are given by
$${\rm Homgr}_A(M,N):=\bigoplus_{i\in \mathbb{Z}}{\rm
Hom}_A(M,N\langle i\rangle),$$ where ${\rm Hom}_A(M,N\langle i
\rangle)$ denotes the space of all graded homomorphisms between
$M$ and $N\langle i\rangle$ (the space of homogeneous morphisms of
degree $i$).

Let $X=(X^i,d^i)$ be a complex of $A$--modules. We say that $X$ is
a complex of graded $A$--modules, or just a graded complex, if for
each $i\in \mathbb{Z}$, $X^i$ is a graded module and $d^i$ is a
homogeneous homomorphism of graded $A$--modules. If $X$ is a
graded complex, then $X\langle j\rangle$ denotes the complex of
graded $A$--modules given by $(X\langle j\rangle)^i:=X^i\langle
j\rangle$ and $d_{X\langle j\rangle}^i:=d^i$. Let $X$ and $Y$ be
graded complexes. A homomorphism $f=\{f^i\}_{i\in\mathbb{Z}}$
between complexes $X$ and $Y$ is a homomorphism of graded
complexes if for each  $i\in \mathbb{Z}$, $f^i$ is a homomorphism
of graded modules.

The category of complexes of graded $A$--modules, denoted by
$C_{gr}(A)$, is the category whose objects are complexes of graded
$A$--modules and morphisms between two graded complexes $X$ and
$Y$ are given by
$${\rm Homgr}_{{\mathcal{C}}_{gr}(A)}(X, Y) := \bigoplus_{i\in\mathbb{Z}}{\rm Homgr}_A(X,Y\langle i\rangle ),$$
where ${\rm Homgr}_A(X,Y\langle i\rangle )$ denotes the space of
graded homomorphisms between $X$ and $Y\langle i\rangle$ (the
space of homogeneous morphisms of degree $i$).

Unless otherwise stated, for a graded algebra $A$, we will assume
that the projective indecomposable $A$-modules are graded in such
way that their tops are in degree 0. For an indecomposable bounded
graded complex of projective $A$-modules, we will assume that the
leftmost non-zero term is graded in such way that its top is in
degree 0.

We say that two symmetric algebras $A$ and $B$ are derived
equivalent if their derived categories of bounded complexes are
equivalent. From Rickard's theory we know that $A$ and $B$ are
derived equivalent if and only if there exists a tilting complex
$T$ of projective $A$--modules such that ${\rm
End}_{K^b(P_A)}(T)\cong B^{op}$.  For more details on derived
categories and derived equivalences we recommend [\ref{KNG}].

\section{Brauer tree algebras}

We now introduce Brauer tree algebras. For a general reference on
Brauer tree algebras we refer reader to [\ref{alp}].

Let $\Gamma$ be a finite connected tree with $e$ edges. We say
that $\Gamma$ is a Brauer tree of type $(m,e)$ if there  is a
cyclic ordering of the edges adjacent to a given vertex, and a
distinguished vertex $v$, called the exceptional vertex, to whom
we assign a positive integer $m$, called the multiplicity of the
exceptional vertex.

Let $A$ be a  finite dimensional symmetric algebra and let
$\Gamma$ be a Brauer tree. We say that  $A$  is a Brauer tree
algebra associated with $\Gamma$ if the isomorphism classes of
simple $A$--modules are in one-to-one correspondence with the
edges of $\Gamma$, and if $P_j$ denotes the projective cover of
the simple $A$--module $S_j$ corresponding to the edge $j$, the
following condition is satisfied:

The heart ${\rm rad }\, P_j/{\rm soc }\, P_j$ is the direct sum of
two uniserial modules $U_v$ and $U_w$, one of which might be zero,
where $v$ and $w$ are vertices of $j$. For $u\in \{v,w\}$, let
$j=j_0, j_1,\dots,j_r$ be the cyclic ordering of the $r+1$ edges
around $u$. The composition factors of $U_u$, starting from the
top, are
$$S_{j_1},S_{j_2},\dots,
S_{j_r},S_{j_0},S_{j_1},S_{j_2},\dots,S_{j_r},S_{j_0},\dots,S_{j_1},S_{j_2},\dots,S_{j_r}$$
where the number of composition factors is $m(r+1)-1$ if $u$ is
the exceptional vertex, and $r$ otherwise.

If $A$ is a Brauer tree algebra associated to a Brauer tree
$\Gamma$ of  type $(m,e)$, we say that $A$ has type $(m,e)$. We
will usually label edges of a Brauer tree by the corresponding
simple modules of a Brauer tree algebra associated to this tree.
We will always assume that the cyclic ordering of the edges
adjacent to a given vertex is given by the counter-clockwise
direction in the given planar embedding of the tree.

\begin{ex}{\rm
A very important Brauer tree of type $(m,e)$ is the Brauer star,
the Brauer tree with $e$ edges adjacent to the exceptional vertex
which has multiplicity $m$.

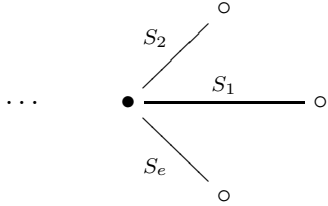
$$\xymatrix{&&\circ&\\
\dots&\bullet\ar@{-}[ur]^{S_2}\ar@{-}[rr]^{S_1}\ar@{-}[dr]_{S_e}&&\circ\\
& &\circ & \\}$$

The composition factors, starting from the top, of  the projective
indecomposable module $P_j$ corresponding to the edge $S_j$ are
$$ S_j, S_{j+1}, \dots, S_e, S_1, S_2,\dots ,S_{j-1}; \,\, \dots \,\, ;S_j, S_{j+1}, \dots, S_e, S_1, S_2,\dots , S_{j-1}; S_j$$
where $S_j$ appears $m+1$ times and $S_i$ appears $m$ times, for
$i\neq j$.

We see that a Brauer tree algebra corresponding to the Brauer star
of type $(m,e)$ is a uniserial algebra. This is very important
because it is easy to do calculations with uniserial modules. }
\end{ex}

Any two Brauer tree algebras associated to the same Brauer tree
$\Gamma$ and defined over the same field $k$ are Morita equivalent
(cf.\ [\ref{Cluj2}], Corollary 4.3.3).

A basic Brauer tree algebra corresponding to a Brauer tree
$\Gamma$ is isomorphic to the algebra $kQ/I$, where $Q$ is a
quiver and $I$ is the ideal of relations. This algebra is
constructed as follows (cf.\ [\ref{Cluj2}], Section 5): We replace
each edge of the Brauer tree by a vertex and for two adjacent
edges, say $j_1$ and $j_2$, which come one after the other in the
circular ordering, say $j_2$ comes after $j_1$, we have an arrow
connecting the two corresponding vertices, starting at the vertex
corresponding to $j_2$ and ending at the vertex corresponding to
$j_1$. If there is only one edge adjacent to the exceptional
vertex and $m>1$, then we add a loop starting and ending at the
vertex corresponding to the edge that is adjacent to the
exceptional vertex. This will give us the quiver $Q$.

Notice that, for each vertex of $\Gamma$ that has more than one
edge adjacent to it,  we get a cycle in the quiver $Q$. The cycle
of $Q$ corresponding to the exceptional vertex will be called the
exceptional cycle. If we assume that $\Gamma$ is not the star,
then the ideal $I$ is generated by two types of relations. The
relations of the first type are given by $ab=0$, where arrows $a$
and $b$ belong to different cycles of $Q$. The second type
relations are the relations of  the form $\alpha=\beta$, for two
cycles $\alpha$ and $\beta$ starting and ending at the same
vertex, neither of which is the exceptional cycle, and  relations
of the form $\alpha^m=\beta$, if $\alpha$ is the exceptional
cycle. The basic algebra $kQ/I$ constructed in this way will be
denoted by $A_{\Gamma}$.

If $\Gamma$ is the star, then the corresponding quiver has only
one cycle and the ideal of relations is generated by all paths of
length $me+1$. This basic algebra corresponding to the Brauer star
will be denoted by $A_S$.

\section{Transfer of gradings}
Let $A$ and $B$ be two symmetric algebras over a field $k$ and let
us assume that $A$ is a graded algebra. The following theorem is
due to Rouquier.
\begin{te}[[\ref{Rou}{]}, {\rm Theorem 6.4}]
Let $A$ and $B$ be as above. Let $T$  be a tilting complex of
$A$-modules that induces derived equivalence between $A$ and $B$.
Then there exists a grading on $B$ and a structure of a graded
complex $T^{\prime}$ on $T$, such  that $T^{\prime}$ induces an
equivalence between the derived categories of graded $A$-modules
and graded $B$-modules.
\end{te}
This theorem tells us that derived equivalences are compatible
with gradings, that is, gradings can be transferred between
symmetric algebras via derived equivalences.

We will now explain how the transfer of gradings via derived
equivalences is done in our context of Brauer tree algebras.
Brauer tree algebras are determined up to derived equivalence by
the multiplicity of the exceptional vertex and the number of edges
of the tree ([\ref{RiStab}], Theorem 4.2). We notice here that the
basic Brauer tree algebra $A_{S}$, which corresponds to the Brauer
star of type $(m,e)$, is naturally graded by putting all arrows in
degree 1. This grading is compatible with radical filtration, in
other words, $A_{S}$ is tightly graded. This means that $A_S$ is
isomorphic to the graded algebra associated with the radical
filtration, i.e.\
$$A_S\cong \bigoplus_{i=0}^{\infty}({\rm rad}\, A_S)^i/({\rm rad}\, A_S)^{i+1}.$$

We will transfer this grading from the algebra $A_S$ to the basic
Brauer tree algebra $A_{\Gamma}$ corresponding to an arbitrary
Brauer tree $\Gamma$ of type $(m,e)$. In order to do that we will
construct a tilting complex of $A_S$-modules which tilts from
$A_S$ to $A_{\Gamma}$. For a given tilting complex $T$ of
$A_S$-modules, which is a bounded complex of finitely generated
projective $A_S$--modules, there exists a structure of a complex
of graded $A_S$--module $T^{\prime}$ on $T$. If $T$ is a tilting
complex that tilts from $A_S$ to $A_{\Gamma}$, then ${\rm
End}_{K^b(P_{A_S})}(T)\cong A_{\Gamma}^{op}$. Viewing $T$ as a
graded complex $T^{\prime}$, and computing this endomorphism ring
as a graded object, we get a graded algebra which is isomorphic to
the opposite algebra of the basic Brauer tree algebra $A_{\Gamma}$
corresponding to $\Gamma$. We notice here that the choice of a
grading on $T^{\prime}$ is unique up to shifting the grading of
each indecomposable summand of $T^{\prime}$.

\subsection{The tilting complex given by Green's walk}

In [\ref{Zak}], the authors give a combinatorial construction of a
tilting complex that tilts from the basic Brauer tree algebra
$A_S$ corresponding to the star of type $(m,e)$, to the basic
Brauer tree algebra $A_{\Gamma}$ corresponding to an arbitrary
Brauer tree of type $(m,e)$. The tilting complexes considered in
[\ref{Zak}] are direct sums of indecomposable complexes which have
no more than two non-zero terms, and a complete classification for
all such tilting complexes is given in [\ref{Zak}].

We will use only one special tilting complex that arises in this
way. It is the complex constructed by taking Green's walk (cf.\
[\ref{Green}]) around $\Gamma$. We construct it as follows:

Starting from the exceptional vertex of a Brauer tree $\Gamma$ of
type $(m,e)$, we take Green's walk around $\Gamma$ and enumerate
all the edges of $\Gamma$. We start this enumeration from an
arbitrary edge adjacent to the exceptional vertex and walk around
$\Gamma$ in the counter-clockwise direction. We will eventually
show that the resulting grading on $A_{\Gamma}$ does not depend on
where we start the enumeration. Define $T$, a tilting complex of
$A_S$-modules, to be the direct sum of the complexes $T_i$, $1\leq
i\leq e$, which correspond to the edges of
 $\Gamma$, and which are defined by induction on the distance of
 an edge from the exceptional vertex in the following way:

\begin{itemize}
\item[(a)] If $i$ is an edge which is adjacent to the exceptional
vertex, then $T_i$ is defined to be the stalk complex $$Q_i\, :\,
0\longrightarrow P_i \longrightarrow 0$$ with $P_i$ in degree 0;

\item[(b)] If $i$ is not adjacent to the exceptional vertex and
assuming  that the shortest path connecting $i$ to the exceptional
vertex is $j_1,j_2,\dots,j_t, i$, where $j_1<j_2<\dots<j_t<i$ in
the labelling from the Green's walk, then $T_i$ is defined to be
the complex $Q_{j_ti}[n_i]$, where $Q_{j_ti}$ is the following
complex with two  non-zero entries in degrees 0 and 1
$$Q_{j_ti}\, : \, 0 \longrightarrow
P_{j_t}\stackrel{h_{j_ti}}{\longrightarrow} P_i\longrightarrow
0,$$ where $h_{j_ti}$ is a homomorphism of the highest possible
rank, and $n_i$ is the shift necessary to ensure that $P_{j_t}$ is
in the same degree as $P_{j_t}$ in other summands of $T$ which are
previously determined.

\end{itemize}

\noindent For the convenience of the reader we include the
following example.
\begin{ex}\label{CompExample}{\rm
Let $\Gamma$ be the following Brauer tree with multiplicity 1 and
with edges numbered by taking Green's walk:
$$
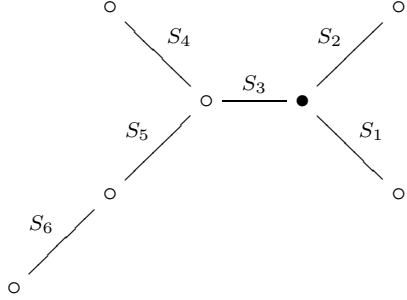
\xymatrix{&\circ\ar@{-}[dr]^{S_4}&&&\circ\\
&&\circ\ar@{-}[r]^{S_3}&\bullet\ar@{-}[dr]^{S_1}\ar@{-}[ur]^{S_2}&\\
&\circ\ar@{-}[ur]^{S_5}&&&\circ\\
\circ\ar@{-}[ur]^{S_6}&&&&\\
}$$

\noindent The tilting complex of $A_S$-modules, where $S$ is the
Brauer star with six edges and multiplicity 1, given by taking
Green's walk around $\Gamma$ is the direct sum
$T=\oplus_{i=1}^6T_i$, where $T_1$, $T_2$ and $T_3$ are the stalk
complexes  $P_1$, $P_2$ and $P_3$ respectively, in degree 0, $T_4$
is $P_3\stackrel{h_{34}}{\longrightarrow} P_4$ with $P_3$ in
degree 0, $T_5$ is $P_3\stackrel{h_{35}}{\longrightarrow} P_5$
with $P_3$ in degree 0, and $T_6$ is
$P_5\stackrel{h_{56}}{\longrightarrow} P_6$ with $P_5$ in degree
1. This complex tilts from $A_S$ to $A_{\Gamma}$, ie.\ ${\rm
End}_{K^b(P_{A_S})}(T)\cong A_{\Gamma}^{op}.$}
\end{ex}

\subsection{Calculating ${\bf {\rm\bf End}_{K^b(P_A)}(T)}$}

Let $A_S$ be a basic Brauer tree algebra corresponding to the star
of type $(m,e)$ and let $A_{\Gamma}$ be a basic Brauer tree
algebra corresponding to a given Brauer tree $\Gamma$ of type
$(m,e)$. Let $T$ be the tilting complex of $A_S$-modules that
tilts from $A_S$ to $A_{\Gamma}$, constructed as in the previous
section, i.e.\ constructed by taking Green's walk around $\Gamma$.
Viewing each summand $T_i$ of $T$ as a graded complex
$T_i^{\prime}$, we have a structure of a graded complex
$T^{\prime}$ on $T$. By calculating ${\rm
Endgr}_{K^b(P_{A_S})}(T^{\prime})\cong A_{\Gamma}^{op}$ we get
$A_{\Gamma}$ as a graded algebra. We will choose $T_i^{\prime}$ to
be $T_i\langle r_i\rangle$, where $r_i$ will be the necessary
shifts that will ensure that the resulting grading is
non-negative. We remind the reader that $T_i$ is assumed to be
graded in such way that its leftmost non-zero term has its top in
degree 0.

We now state the main theorem of this section.
\begin{te}\label{glteorema}
Let $\Gamma$ be an arbitrary Brauer tree with $e$ edges and
multiplicity $m$ of the exceptional vertex and  let $A_{\Gamma}$
be the basic Brauer tree algebra determined by this tree. The
algebra $A_{\Gamma}$ can be non-negatively graded.
\end{te}
\noindent {\bf Proof.} In order to grade $A_{\Gamma}$, we need to
calculate ${\rm Homgr}_{K^b(P_{A_S})}(T_i^{\prime},T_j^{\prime})$,
as   graded vector spaces, for those $T_i^{\prime}$ and
$T_j^{\prime}$ which correspond to edges $S_i$ and $S_j$ that are
adjacent to the same vertex of $\Gamma$, and which come one after
the other,  $i$ after $j$, in the circular ordering associated to
that vertex. This is a consequence of the fact that when
identifying ${\rm End}_{K^b(P_{A_S})}(T)^{op}$ with $A_{\Gamma}$,
elements corresponding to the vertices of the quiver of
$A_{\Gamma} $ are given by  $ {\rm id}_{T_i}\in  {\rm
End}_{K^b(P_{A_S})}(T_i)$, $i=1,2,\dots,e $; and the subspace of
$A_{\Gamma}$ generated by all paths starting at the vertex of
$A_{\Gamma}$ corresponding  to ${\rm id}_{T_i}\in {\rm
End}_{K^b(P_{A_S})}(T_i)$ and finishing at the vertex
corresponding to ${\rm id}_{T_j} \in {\rm
End}_{K^b(P_{A_S})}(T_j)$, is given by ${\rm
Hom}_{K^b(P_{A_S})}(T_i, T_j)$. In fact, we only need to calculate
the non-zero summand of ${\rm
Homgr}_{K^b(P_{A_S})}(T_i^{\prime},T_j^{\prime})$ which is in the
lowest degree. That will be degree of the unique arrow of the
quiver of $A_{\Gamma}$ pointing from the vertex corresponding to
$S_i$ to the vertex corresponding to $S_j$.

The dimension of  ${\rm Hom}_{K^b(A_S)}(T_i,T_j)$
\begin{itemize} \item[{\rm (1)}] is $0$, if the vertices
corresponding to ${\rm id}_{T_i}$ and ${\rm id}_{T_j}$ do not
belong to the same cycle, \item[{\rm (2)}]
\begin{itemize}\item[{\rm (2.a)}]is $m$, if $i\neq j$ and the
vertices corresponding to ${\rm id}_{T_i}$ and ${\rm id}_{T_j}$
belong to the exceptional cycle, \item[{\rm (2.b)}] is $1$, if
 $i\neq j$ and the vertices corresponding to ${\rm id}_{T_i}$ and
${\rm id}_{T_j}$ belong to the same non-exceptional cycle,
\end{itemize} \item[{\rm (3)}] \begin{itemize}\item[{\rm (3.a)}]
is $m+1$, if $i=j$ and the vertex corresponding to ${\rm
id}_{T_i}$ belongs to the exceptional cycle, \item[{\rm (3.b)}]is
$2$,  if $i=j$ and the vertex corresponding to ${\rm id}_{T_i}$
does not belong to the exceptional cycle.
\end{itemize}
\end{itemize}

Edges $i$ and $j$ of a Brauer tree $\Gamma$ that are adjacent to
the same vertex, say $v$, and that come one after the other in the
circular ordering of $v$, can either be at the same distance from
the exceptional vertex or the distance of one of them, say $i$, is
one less than the distance of $j$. If the former holds, then a
part of $\Gamma$ is given by
$$
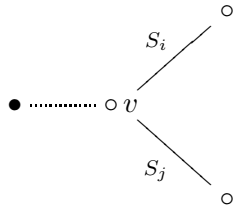
\xymatrix{
&&&\circ\\&\bullet\ar@{.}[r]&\circ \,
v\ar@{-}[ur]^{S_i}\ar@{-}[dr]_{S_j}&\\&&&\circ}$$ where $v$ may or
may not be the exceptional vertex.

If the latter holds, than a part of $\Gamma$ is given by one of
the following two diagrams
$$
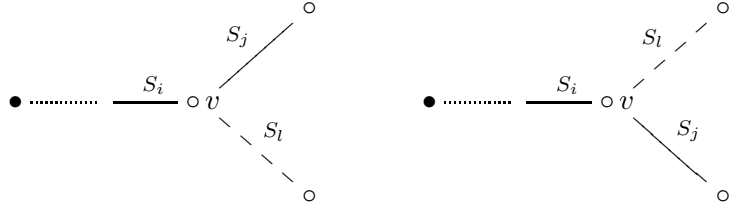
\xymatrix{&&&\circ\\
\bullet\ar@{.}[r]&\ar@{-}[r]^{S_i}& \circ\, v \ar@{-}[ur]^{S_j}\ar@{--}[dr]^{S_l}&\\
&&&\circ }\quad \quad \quad \xymatrix{&&&\circ\\
\bullet\ar@{.}[r]&\ar@{-}[r]^{S_i}& \circ\, v \ar@{--}[ur]^{S_l}\ar@{-}[dr]^{S_j}&\\
&&&\circ }
$$
where the leftmost vertex of $S_i$ may or may not be the
exceptional vertex, and there may or may not be more edges, such
as $S_l$, adjacent to $v$. In the case of the first diagram we
have an arrow from $S_i$ to $S_j$ in the quiver of $A_{\Gamma}$,
and in the case of the second diagram we have an arrow from $S_j$
to $S_i$.

It follows that it is sufficient to consider the following four
cases:

{\textbf{CASE 1}}   Edges $S_i$ and $S_j$ are adjacent to the
exceptional vertex. The corresponding part of $\Gamma$ is
$$
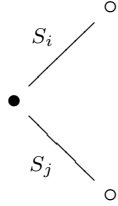
\xymatrix{
&&\circ\\&\bullet\ar@{-}[ur]^{S_i}\ar@{-}[dr]_{S_j}&\\&&\circ}
$$
In this case, the corresponding summands of the graded tilting
complex $T^{\prime}$  are $T_i^{\prime}:=T_i$ and
$T_j^{\prime}:=T_j$, where $T_i$ and $T_j$ are stalk complexes
$P_i$ and $P_j$ concentrated in degree 0. If $i>j$, then  ${\rm
Homgr}_{K^b(P_A)}(T_i^{\prime},T_j^{\prime})\cong {\rm
Homgr}_{A_S}(P_i, P_j)\cong k\langle -(i-j)\rangle \oplus M$, as
graded vector spaces, where $M$ is the sum of the summands that
appear in higher degrees than $i-j$. In other words, $i-j$ is the
degree of the corresponding arrow of the quiver of $A_{\Gamma}$
whose source is $S_i$ and whose target is $S_j$. If $i<j$, then
the degree of the corresponding arrow of the quiver of
$A_{\Gamma}$ is $e-(j-i)$, where $e$ is the number of edges of the
tree.

If there is only one edge adjacent to the exceptional vertex and
$m>1$, then the corresponding loop will be in degree $e$.

\textbf{CASE 2} Edges $S_i$ and $S_j$ are adjacent to a
non-exceptional vertex $v$ and one of them, say $S_i$, is adjacent
to the exceptional vertex. Then we have that a part of $\Gamma$ is
given by one of the following two diagrams
$$
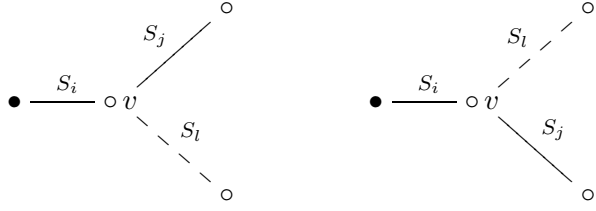
\xymatrix{&&\circ\\
\bullet\ar@{-}[r]^{S_i}="a"& \circ\, v \ar@{-}[ur]^{S_j}="b"\ar@{--}[dr]^{S_l}&\\ 
&&\circ }\quad\quad\quad\quad \xymatrix{&&\circ\\
\bullet\ar@{-}[r]^{S_i}& \circ\, v \ar@{--}[ur]^{S_l}\ar@{-}[dr]^{S_j}&\\
&&\circ }
$$
where it could happen that there are more edges adjacent to $v$,
such as $S_l$.

In this case, the summands of the tilting complex $T^{\prime}$ are
$T_i^{\prime}:=Q_i$, the stalk complex with $P_i$ in degree 0, and
$T_j^{\prime}:=Q_{ij}$, where $Q_{ij}$ is previously defined
complex $P_i\stackrel{h_{ij}}{\longrightarrow} P_j$, with $P_i$ in
degree 0.

Since $i<j$, we have that any map from $Q_i$ to $Q_{ij}$ has to
map $P_i$ to the kernel of the map $h_{ij}$. It follows that this
map has to map ${\rm top}\,P_i$ to ${\rm soc} \,P_i$. This means
that the corresponding arrow of the quiver of $A_{\Gamma}$, whose
source is $S_i$ and whose target is $S_j$, is in degree $me$. This
happens in case of the first diagram.

In case of the second diagram,  there is an arrow from $S_j$ to
$S_i$ in the quiver of $A_{\Gamma}$. The identity map on $P_i$
will give us a morphism between graded complexes $T_j^{\prime}$
and $T_i^{\prime}$ and we conclude that the corresponding arrow
from $S_j$ to $S_i$ is in degree 0.

\textbf{CASE 3}  Edges $S_i$ and $S_j$, where $i<j$, are adjacent
to a non-exceptional vertex $v$, $S_i$ is closer to the
exceptional vertex than $S_j$, and the leftmost vertex of $S_i$ is
non-exceptional. Then we have that a part of $\Gamma$ is given by
one of the following two diagrams
$$
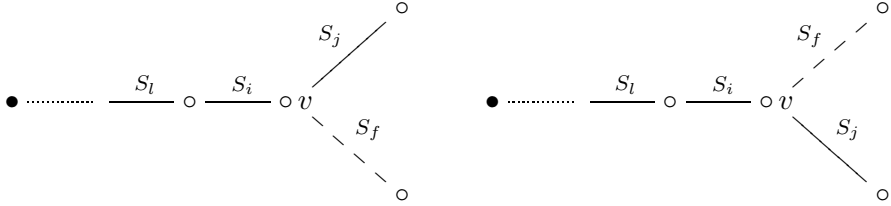
\xymatrix{&&&&\circ\\
\bullet\ar@{.}[r]&\ar@{-}[r]^{S_l}&\circ \ar@{-}[r]^{S_i}& \circ\, v \ar@{-}[ur]^{S_j}\ar@{--}[dr]^{S_f}&\\
&&&&\circ } \quad\quad \xymatrix{&&&&\circ\\
\bullet\ar@{.}[r]&\ar@{-}[r]^{S_l}&\circ \ar@{-}[r]^{S_i}& \circ\, v \ar@{--}[ur]^{S_f}\ar@{-}[dr]^{S_j}&\\
&&&&\circ }
$$
where it could happen that there are more edges adjacent to $v$,
such as $S_f$. The leftmost vertex of $S_l$ can be either
exceptional or non-exceptional.

The summands $T_i^{\prime}$ and $T_j^{\prime}$ of the graded
tilting complex $T^{\prime}$ corresponding to the edges $S_i$ and
$S_j$ are defined to be graded complexes $(Q_{li}[-r_{li}])\langle
r_{li}\rangle$ and $(Q_{ij}[-r_{li}-1])\langle r_{li}+1 \rangle$,
where we set $r_{li}$ to be  the distance between the exceptional
vertex and the leftmost vertex of $S_l$. The horizontal shifts
$[-r_{li}]$ and $[-r_{li}-1]$ are necessary to ensure that $P_i$
appears in the same degree as $P_i$ in all other previously
defined summands of $T$. The vertical shifts $\langle
r_{li}\rangle$ and $\langle r_{li}+1\rangle$ are necessary to
ensure that top of the rightmost term of $Q_{li}$, which is $P_i$,
and top of the leftmost term of $Q_{ij}$, which is also $P_i$, are
in the same degree. This way we avoid negative degrees.

If the first diagram occurs, then any morphism of graded complexes
from $(Q_{li}[-r_{li}])\langle r_{li} \rangle$ to
$(Q_{ij}[-r_{r_{li}-1}])\langle r_{li}+1 \rangle$ has to map ${\rm
top}\, P_i$ to ${\rm soc}\, P_i$. From this we conclude that the
corresponding arrow of the quiver of $A_{\Gamma}$ that points from
$S_i$ to $S_j$ is in degree $me$.

If the second diagram occurs, then the identity map on $P_i$ gives
us a map from $(Q_{ij}[-r_{li}-1])\langle r_{li}+1 \rangle$ to
$(Q_{li}[-r_{li}])\langle r_{li} \rangle$ which is not homotopic
to zero. From this we have that the arrow from $S_j$ to $S_i$ is
in degree 0.

\textbf{CASE 4}   Edges $S_i$ and $S_j$, where $j<i$, are adjacent
to a non-exceptional vertex $v$, and the edge with minimal index
among the edges adjacent to $v$, say $S_l$, comes  before $S_j$ in
the circular ordering of $v$. Then we have that a part of $\Gamma$
is
$$
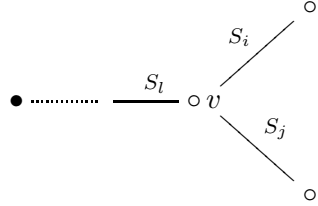
\xymatrix{&&&\circ\\
\bullet\ar@{.}[r]&\ar@{-}[r]^{S_l}&\circ\, v\ar@{-}[ur]^{S_i}\ar@{-}[dr]^{S_j}&\\
&&&\circ }
$$
\noindent where the leftmost vertex of $S_l$ can be either
exceptional or non-exceptional. In each case, the summands
$T_i^{\prime}$ and $T_j^{\prime}$ of the graded tilting complex
$T^{\prime}$ corresponding to the edges $S_i$ and $S_j$ are
$(Q_{li}[-r_{li}])\langle r_{li} \rangle$ and
$(Q_{lj}[-r_{li}])\langle r_{li} \rangle$, where we again set
$r_{li}$ to be the distance between the exceptional vertex and the
leftmost vertex of $S_l$. The map $({\rm id}_{P_l}; h_{ij})$,
where $h_{ij}$ is a map of the maximal rank,  from
$(Q_{li}[-r_{li}])\langle r_{li} \rangle$ to
$(Q_{lj}[-r_{li}])\langle r_{li} \rangle$ will give us a nonzero
map in ${\rm Homgr}_{K^b(P_A)}(T_i^{\prime},T_j^{\prime})$ and
this map is in degree 0. Therefore, the corresponding arrow from
$S_i$ to $S_j$ is in degree 0. \vspace{2mm}

These four cases cover all the possible local structures of a
Brauer tree $\Gamma$ that we can encounter when putting a grading
on the basic Brauer tree algebra $A_{\Gamma}$ that corresponds to
this tree. We only need to walk around the Brauer tree $\Gamma$
and recognize which of the four cases occurs for the adjacent
edges $S_i$ and $S_j$. In each of the four cases above, we have
that the corresponding arrows are in non-negative degrees. Hence,
the resulting grading on $A_{\Gamma}$ is non-negative.
$\blacksquare$

The grading on $A_{\Gamma}$ constructed in the proof of the
previous theorem will be referred to as the grading constructed by
taking Green's walk around $\Gamma$.

\begin{ex}\label{linePR}{\rm  (a) Let $\Gamma$ be the Brauer
tree from  Example \ref{CompExample}. We first construct the
quiver of the basic Brauer tree algebra $A_{\Gamma}$ corresponding
to this tree. Each edge is replaced by a vertex and for two
adjacent edges, which come one after the other in the circular
ordering, we have an arrow connecting  two corresponding vertices
of the quiver in the opposite ordering of the circular ordering of
edges, i.e.\ in the clockwise direction.
$$
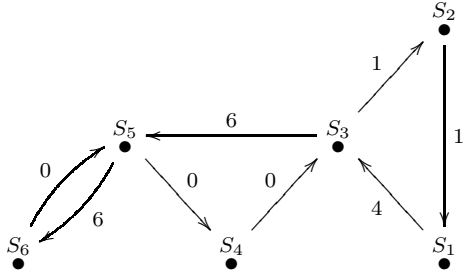
\xymatrix{
&&&&\stackrel{S_2}{\bullet}\ar[dd]^{1}\\
&\stackrel{S_5}{\bullet}\ar@/^/[dl]^{6}\ar[dr]^{0}&&\stackrel{S_3}{\bullet}\ar[ll]_{6}\ar[ur]^{1}&\\
\stackrel{S_6}{\bullet}\ar@/^/[ur]^{0}&&\stackrel{S_4}{\bullet}\ar[ur]^{0}&&\stackrel{S_1}{\bullet}\ar[ul]^{4}\\
}
$$
The degrees of the arrows between $S_1$ and $S_3$, $S_3$ and
$S_2$, $S_2$ and $S_1$ are computed using the first case from the
proof of the previous theorem. The degree of the arrow between
$S_3$ and $S_5$ is 6 and is computed using the second case. Arrows
between $S_5$ and $S_4$, and $S_4$ and $S_3$ are in degree 0 and
these degrees are computed as in the fourth and the second case
respectively. The degree of the arrow between $S_5$ and $S_6$ is 6
and the degree of the arrow between $S_6$ and $S_5$ is 0 and these
are computed as in the third case.

(b) If the Brauer tree is the line with $e$ edges and the
exceptional vertex at the end, ie.\

$$
\xymatrix{\bullet\ar@{-}[r]^{S_1}&\circ\ar@{-}[r]^{S_2}&\circ\ar@{-}[r]^{S_3}&\circ\ar@{-}[r]^{S_{4}}&\dots\ar@{-}[r]^{S_{n-1}}&\circ\ar@{-}[r]^{S_e}&\circ}
$$
then the basic Brauer tree algebra $A_{\Gamma}$ is graded and its
quiver has $e$ vertices
$$
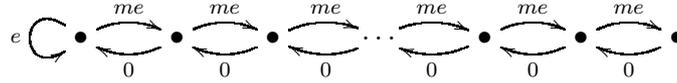
\xymatrix{\bullet\ar@(ul,dl)_{e}\ar@/^/[r]^{me}&\bullet\ar@/^/[r]^{me}\ar@/^/[l]^0&\bullet\ar@/^/[l]^0\ar@/^/[r]^{me}&\dots\ar@/^/[l]^0\ar@/^/[r]^{me}&\bullet\ar@/^/[l]^0\ar@/^/[r]^{me}&\bullet\ar@/^/[l]^0\ar@/^/[r]^{me}&\bullet\ar@/^/[l]^0}
$$
\noindent If $m=1$, then there is no loop in the above quiver.}
\end{ex}
\vspace{2mm} Let $\Gamma$ be an arbitrary Brauer tree. Each edge
of $\Gamma$ that is adjacent to the exceptional vertex determines
a connected subtree of a Brauer tree $\Gamma$. We  call these
subtrees components of the Brauer tree.
\begin{lm}\label{NoEdgesComp}
Let $\alpha$ be an arrow contained in the exceptional cycle of the
quiver of $A_{\Gamma}$ which starts at $S_i$ and ends at $S_j$. If
$A_{\Gamma}$ is graded by taking Green's walk around $\Gamma$,
then the degree of $\alpha$ is equal to the number of edges in the
component of  $\Gamma$ corresponding to $S_j$.
\end{lm}
\noindent {\bf Proof.} If $i>j$, then because of the way we
enumerate edges by taking Green's walk we have that $i=j+s$, where
$s$ is the number of edges in the component corresponding to
$S_j$. Hence,  $\alpha$ is in degree $i-j=s$.

If $i<j$, which only happens if $i=1$, then $\alpha$ is in degree
$e-(j-i)=e-j+1$, and this number is equal to the number of edges
of the component corresponding to $S_j$. $\blacksquare$

From this lemma it follows that the resulting grading of the
exceptional cycle does not depend on from which edge adjacent to
the exceptional vertex we start enumeration of edges. This leads
us to the following proposition.

\begin{prop}
Let $A_\Gamma$ be the basic Brauer tree algebra whose tree is
$\Gamma$ and let us assume that $A_\Gamma$ is graded by taking
Green's walk around $\Gamma$. The resulting grading does not
depend on from which edge adjacent to the exceptional vertex we
start Green's walk.
\end{prop}
\noindent {\bf Proof.} Let us assume that we have done two walks
around $\Gamma$ starting at a different edge each time. Let us
assume that the index of  $S$, where $S$ is one of the edges
adjacent to the exceptional vertex, is 1 in the first walk, and
that its index is $1+l$ in the second walk.  Let us assume that we
got two tilting complexes $T_1$ and $T_2$ of $A_S$-modules by
taking these two walks. These complexes are equal up to cyclic
permutation of the vertices of the Brauer star $A_S$. In other
words, each index of each term (which is a projective
indecomposable $A_S$-module) of each summand of $T_1$ has been
cyclically permuted by $l$ to get the corresponding index of the
corresponding term of the corresponding summand of $T_2$. These
two tilting complexes will give us the same grading because of the
'cyclic' structure of the Brauer star $A_S$. $\blacksquare$

\begin{lm}\label{lema1} Let $A_\Gamma$ be the basic Brauer tree algebra whose
tree is $\Gamma$. Let $Q$ be its quiver and let us assume that
$A_\Gamma$ is graded by taking Green's walk around $\Gamma$.  The
only cycle of $Q$ that does not contain any arrows that are in
degree $0$ is the exceptional cycle. For a non-exceptional cycle
there is exactly one arrow that is not in degree $0$. This arrow
is in degree $me$ and the index of its target is greater than the
index of its source.
\end{lm}
\noindent{\bf Proof.} If $\alpha$ is an arbitrary arrow of the
exceptional vertex, then $\alpha$ is in a positive degree. This
follows from Lemma \ref{NoEdgesComp}, since the number of edges in
each component is strictly positive. For a non-exceptional cycle,
from the last three cases from the proof of Theorem
\ref{glteorema}, we see that exactly one arrow is not in degree
$0$. This is the arrow whose source is the vertex of that cycle
with the minimal index, and whose target is the vertex of that
cycle with the maximal index. Its degree is $me$ by the proof of
Theorem \ref{glteorema}. $\blacksquare$

For the arrows of the quiver of $A_{\Gamma}$ that are in degree 0,
we have that the index of their source is greater than the index
of their target. We state this in the following lemma.

\begin{lm}\label{indicesArr}
Let $\alpha$ be an arrow of a non-exceptional cycle which is in
degree $0$. Then the index of the source of $\alpha$ is greater
than the index of the target of $\alpha$.
\end{lm}

\begin{lm}
Let $A_{\Gamma}$, $\Gamma$ and $Q$ be as above. The socle of
$A_{\Gamma}$ is in degree $me$.
\end{lm}
\noindent {\bf Proof.} For an arbitrary cycle, say $\gamma$, let
$\alpha_1,\dots,\alpha_r$ be the arrows of that cycle, appearing
in that cyclic ordering. If $\gamma$ is a non-exceptional cycle of
the quiver $Q$, then paths of the form
$\alpha_i\alpha_{i+1}\dots\alpha_{i+r-1}$, where the addition in
indices is modulo $r$, and $1\leq i\leq r$, belong to the socle.
If $\gamma$ is  the exceptional cycle, then the paths
$(\alpha_i\alpha_{i+1}\dots\alpha_{i+r-1})^m$ belong to the socle.
These elements span the whole socle. For a non-exceptional cycle,
the only arrow which is in a non-zero degree is the arrow whose
source is the vertex of that cycle with the minimal index, and
whose target is the vertex of that cycle with the maximal index.
By the cases 2,3,4 from the proof of Theorem \ref{glteorema}, that
arrow is in degree $me$. It follows that the path
$\alpha_i\alpha_{i+1}\dots\alpha_{i+r-1}$, for all $i$, is in
degree $me$. For the exceptional cycle, let us assume that the
vertices contained in the exceptional cycle are $S_{l_1}, S_{l_2},
\dots, S_{l_r}$, in that cyclic ordering, and that
$l_1>l_2>\dots>l_r=1$. The sum of degrees of the arrows of the
exceptional cycle is
$\sum_{j=1}^{r-1}(l_{j}-l_{j+1}))+(e-(l_1-l_r))=e$. (This also
follows from Lemma \ref{NoEdgesComp}) Therefore, the $m^{th}$
power of $\alpha_i\alpha_{i+1}\dots\alpha_{i+r-1}$, for all $i$,
is in degree $me$. $\blacksquare$

\section{The subalgebra $A_0$}
Let  $\Gamma$ be a given Brauer tree and let $A_{\Gamma}$ be a
basic Brauer tree algebra associated with this tree.  Let $T$ be
the tilting complex that we constructed by taking Green's walk
around $\Gamma$. If we assume that we have a graded algebra
$A_{\Gamma}$ by using this complex, then the subalgebra $A_0$
consisting of the elements that are in degree 0 has an interesting
structure.

The quiver of the basic algebra $A_{\Gamma}$ is the union of the
cycles contained in it. From Lemma \ref{lema1} it follows that the
only cycle that does not contain any arrows that are in degree 0
is the exceptional cycle. If we assume that there are $t$ edges
that are adjacent to the exceptional vertex of $\Gamma$, we see
immediately that the exceptional cycle divides the quiver of the
subalgebra $A_0$ into $t$ disjoint parts, each labelled by a
vertex of the exceptional cycle corresponding to an edge adjacent
to the exceptional vertex of $\Gamma$.

\begin{prop}
Let $A_{\Gamma}$  be a basic Brauer tree algebra associated with a
given Brauer tree $\Gamma$ and let $T$, $A_0$ and $t$ be as above.
The algebra $A_0$ is the direct product of $t$ subalgebras.
\end{prop}
\noindent{\bf Proof.} The factors of $A_0$ are path algebras of
$t$ disjoint subquivers of the quiver of $A_0$. $\blacksquare$

Each of these factors in the previous proposition is labelled by
the corresponding vertex which belongs to the exceptional cycle.
Let $A_v$ be the connected component of $A_0$ that corresponds to
a vertex $v$ of the exceptional cycle.

\begin{lm}
In the quiver of the component $A_v$ of $A_0$ there is at most one
arrow with vertex $v$ as its target. For any other vertex of the
quiver of $A_v$ there are at most two arrows with that vertex as a
target.
\end{lm}
\noindent {\bf Proof.} In the quiver $Q$ of the basic Brauer tree
algebra $A_{\Gamma}$, each vertex is contained in at most two
cycles. Hence, for an arbitrary vertex $v$ of $Q$, there are at
most two arrows of $Q$ whose target is $v$.  If one of the cycles
that contain $v$ is the exceptional cycle, then one of those two
arrows whose target is $v$ is in a positive degree. Therefore, for
a given vertex $v$ of the exceptional cycle, there is at most one
arrow which is in degree 0 that has $v$ as its target. Also, for
every other vertex in this component of the quiver of $A_0$, there
are at most two arrows with that vertex as a target. From Lemma
\ref{indicesArr}, we see that these arrows point from a vertex
with a larger index to a vertex with a smaller index.
$\blacksquare$

\begin{lm}
Let $w$ be a vertex of the quiver of $A_v$ different from $v$. In
the quiver of $A_v$ there is exactly one arrow that has $w$ as its
source.
\end{lm}
\noindent {\bf Proof.} The vertex $w$ belongs either to one or to
two cycles of $Q$, depending on whether the corresponding edge of
the Brauer tree is an end edge or not. Therefore, there are either
one or two arrows that have $w$ as its source. If there are two
such cycles, then the corresponding edge is not an end edge. Then
the arrow that has $w$ as its source and that has vertex of a
greater index than $w$ as its target is in a positive degree by
Lemma \ref{indicesArr}. The other arrow of $Q$ that has $w$ as its
source and a vertex of a smaller index than the index of $w$ as
its target is in degree 0, by the same lemma. $\blacksquare$

\begin{prop}
The quiver of $A_v$ is a directed rooted tree with vertex $v$ as
its root, and with arrows pointing from higher levels of the tree
to  lower levels of the tree, with root $v$ being in level $0$.
\end{prop}
\noindent {\bf Proof.} From the previous two lemmas we conclude
easily that the component of the quiver of $A_0$ does not have any
cycles, because one arrow in each non-exceptional cycle of the
quiver of $A_{\Gamma}$ is in a positive degree. Also, it follows
that this component is a tree with at most two arrows having the
same target, and at most one arrow having an arbitrary vertex as
its source. If we view this tree as a rooted tree with the vertex
that belongs to the exceptional cycle as the root, then all arrows
point from the higher levels to the lower levels by Lemma
\ref{indicesArr}, with the root being in level 0. $\blacksquare$

\begin{prop}
Each of the components of the subalgebra $A_0$ is the path algebra
of a directed rooted tree with arrows pointing from  higher levels
towards lower levels.  The only relations that occur in these
components are of the form $\alpha\beta=0$, where $\alpha$ and
$\beta$ are arrows that belong to different cycles of the quiver
of $A_{\Gamma}$, such that the target of $\alpha$ is the source of
$\beta$.
\end{prop}
\noindent {\bf Proof.} It is left to prove that we only have
relations of type $\alpha\beta=0$. These relations are inherited
from the relations of the algebra $A_{\Gamma}$. The only other
relations that appear in  $A_{\Gamma}$ are of type $\rho=\sigma$,
where $\rho$ and $\sigma$ are two cycles having the same source
and target. Since in the quiver of  $A_0$ there are no cycles,
these relations are not present. $\blacksquare$
\begin{co}
The subalgebra $A_0$ is tightly graded.
\end{co}
\noindent {\bf Proof.} By the previous proposition, the ideal of
relations of $A_0$ is generated by elements of the form
$\alpha\beta$. If the arrows of the quiver of $A_0$ are in degree
1, then these generators are homogeneous of degree 2. The ideal of
relations is homogeneous, hence the quotient algebra $A_0$ of the
path algebra $kQ$ is also graded with arrows in degree 1, i.e.\ it
is tightly graded. $\blacksquare$

\begin{ex}\label{NulaPr}
{\rm If $\Gamma$ is the following Brauer tree with multiplicity 1
and edges numbered by taking Green's walk:

$$
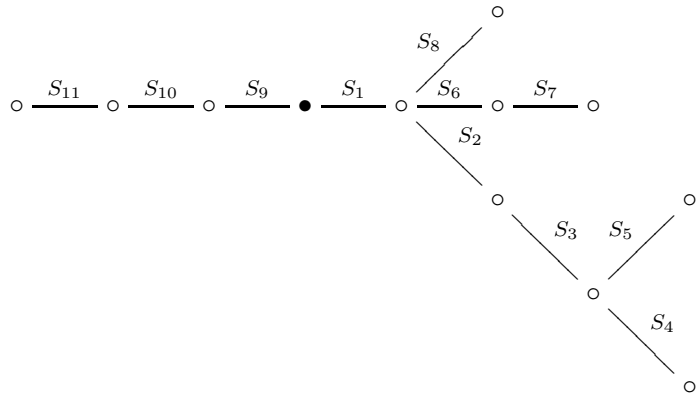
\xymatrix{&&&&&\circ &&\\
\circ \ar@{-}[r]^{S_{11}}&\circ\ar@{-}[r]^{S_{10}}&\circ\ar@{-}[r]^{S_9}&\bullet\ar@{-}[r]^{S_1}&\circ\ar@{-}[dr]^{S_2}\ar@{-}[ur]^{S_8}\ar@{-}[r]^{S_6}&\circ\ar@{-}[r]^{S_7}&\circ &\\
&&&&&\circ\ar@{-}[dr]^{S_3}&&\circ\\
&&&&&&\circ\ar@{-}[dr]^{S_4}\ar@{-}[ur]^{S_5}&\\
&&&&&&&\circ}
$$
then the quiver of the graded basic Brauer tree algebra
$A_{\Gamma}$ associated with this tree is

$$
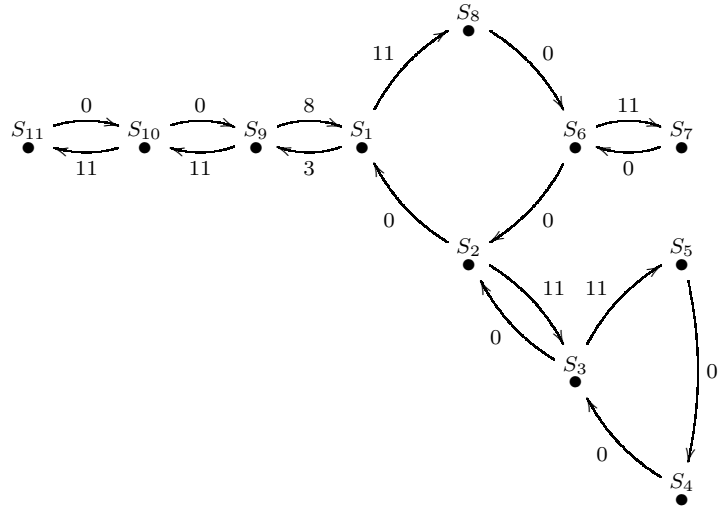
\xymatrix{
&&&&\stackrel{S_8}{\bullet}\ar@/^/[dr]^0&&\\
\stackrel{S_{11}}{\bullet}\ar@/^/[r]^0&\stackrel{S_{10}}{\bullet}\ar@/^/[l]^{11}\ar@/^/[r]^0&\stackrel{S_9}{\bullet}\ar@/^/[l]^{11}\ar@/^/[r]^8&\stackrel{S_1}{\bullet}\ar@/^/[l]^3\ar@/^/[ur]^{11}&&\stackrel{S_6}{\bullet}\ar@/^/[dl]^0\ar@/^/[r]^{11}&\stackrel{S_7}{\bullet}\ar@/^/[l]^0\\
&&&&\stackrel{S_2}{\bullet}\ar@/^/[dr]^{11}\ar@/^/[ul]^0&&\stackrel{S_5}{\bullet}\ar@/^/[dd]^0\\
&&&&&\stackrel{S_3}{\bullet}\ar@/^/[ul]^0\ar@/^/[ur]^{11}&\\
&&&&&&\stackrel{S_4}{\bullet}\ar@/^/[ul]^0 }
$$
The algebra $A_0$ is consisted of two components because there are
two edges adjacent to the exceptional vertex.  The quiver of the
first component is
$$\xymatrix{\bullet\ar[r]^{a_1}&\bullet\ar[r]^{a_2}&\bullet}$$ and
the only relation is $a_1a_2=0.$ The quiver of the second
component is
$$
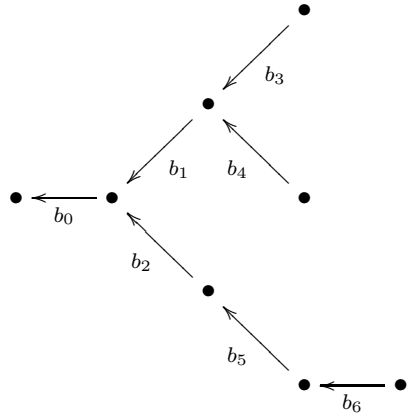
\xymatrix{
&&&\bullet\ar[dl]^{b_3}&\\
&&\bullet\ar[dl]^{b_1}&&\\
\bullet&\bullet\ar[l]^{b_0}&&\bullet\ar[ul]^{b_4}&\\
&&\bullet\ar[ul]^{b_2}&&\\
&&&\bullet\ar[ul]^{b_5}&\bullet\ar[l]^{b_6}}
$$
and the relations are $b_2b_0=0$, $b_5b_2=0$ and $b_4b_1=0.$ }
\end{ex}

\subsection{Recovering the quiver of $A_{\Gamma}$ from the quiver of $A_0$}
The grading resulting from taking Green's walk has some
interesting properties. We will see in this section that the
algebra $A_0$ carries a lot of information about the algebra
$A_{\Gamma}$.

Let $\Gamma$, $A_{\Gamma}$ and $A_0$ be as before. If we omit the
arrows of the exceptional cycle of the quiver of $A_{\Gamma}$, we
see that the resulting quiver  consists of the connected
components which correspond to the connected components of the
quiver of $A_0$.  If we look at the components of $A_0$ we see
that it is sufficient to know the quiver and relations of such
component to recover the quiver of the corresponding component of
$A_{\Gamma}$. This is a consequence of Lemma \ref{lema1}, which
says that in every non-exceptional cycle of the quiver of
$A_{\Gamma}$ there is only one arrow that is not in degree 0.

Let $Q_v$ be one of the connected components of the quiver of
$A_0$ and let $Q_1$ be the corresponding component of the quiver
of $A_{\Gamma}$. We have seen in the previous section that $Q_v$
is a rooted tree with the root $v$ belonging to the exceptional
cycle. Starting from the root of this tree, we can recover the
corresponding component $Q_1$ of the quiver $Q$.

\begin{prop}
Let $Q_v$  be a connected component of the quiver of  $A_0$ and
let $Q_1$ be its corresponding connected component of the quiver of
$A_{\Gamma}$ when the exceptional cycle is omitted. From the
quiver $Q_v$ and its relations we can recover the quiver $Q_1$ and
the relations of the corresponding component of the algebra
$A_{\Gamma}$.
\end{prop}
\noindent {\bf Proof.} Start from the root $v$ of the rooted tree
$Q_v$. Take the longest non-zero path, say $\rho$, ending at $v$.
Add an arrow pointing from $v$ to the source vertex of $\rho$. If
there is no such path of length greater than 1, then add an arrow
from $v$ to the starting point of the only arrow ending at the
root $v$. In this way we recover the cycle of $Q_1$ which has root
$v$ as one of its vertices. The added arrow was an arrow of $Q_1$
that is in a non-zero degree. Now, we repeat the same step with an
arbitrary vertex in the level 1 of the rooted tree instead of the
root, but we only consider paths which do not contain arrows that
belong to already recovered cycles. Repeat the same step for all
other vertices in level 1 of the rooted tree $Q_v$. Repeat the
same steps for vertices in other levels of the rooted tree $Q_v$
until every cycle is recovered.  In this way we recovered the
whole corresponding component $Q_1$  of the quiver of
$A_{\Gamma}$.  As far as the relations are concerned, we get
relations for the basic Brauer tree algebra corresponding to a
given tree, ie.\ for two successive arrows belonging to two
different cycles we set their product to be zero and we set two
cycles starting and ending at the same vertex to be equal.
$\blacksquare$

\begin{ex}
 {\rm
Let $A_{\Gamma}$ be a basic Brauer tree algebra corresponding to
the Brauer tree from  Example \ref{NulaPr}. The algebra $A_0$ has
two components. Let us recover the corresponding components of
$A_{\Gamma}$. The first component is given by the quiver
$$\xymatrix{\bullet\ar[r]^{a_1}&\bullet\ar[r]^{a_2}&\bullet}$$ and the relation $a_1a_2=0$.
Starting from the root we immediately recover the first cycle
since there is no non-zero path of length greater than 1 whose
target is the root. Consequently, the second cycle is easily
recovered and we get that the corresponding component of the
quiver of $A_{\Gamma}$ is
$$
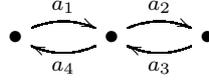
\xymatrix{\bullet\ar@/^/[r]^{a_1}&\bullet\ar@/^/[l]^{a_4}\ar@/^/[r]^{a_2}&\bullet\ar@/^/[l]^{a_3}}$$

The quiver of the second component is
$$
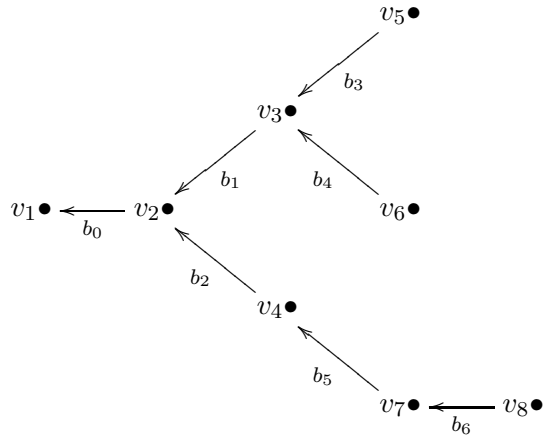
\xymatrix{
&&&v_5\bullet\ar[dl]^{b_3}&\\
&&v_3\bullet\ar[dl]^{b_1}&&\\
v_1\bullet&v_2\bullet\ar[l]^{b_0}&&v_6\bullet\ar[ul]^{b_4}&\\
&&v_4\bullet\ar[ul]^{b_2}&&\\
&&&v_7\bullet\ar[ul]^{b_5}&v_8\bullet\ar[l]^{b_6}}
$$
and the relations are $b_2b_0=0$, $b_5b_2=0$ and $b_4b_1=0.$ The
longest non-zero path ending at $v_1$ is $b_3b_1b_0$. Therefore we
have to add an arrow from $v_1$ to $v_5$. This will give us the
following partial quiver
$$
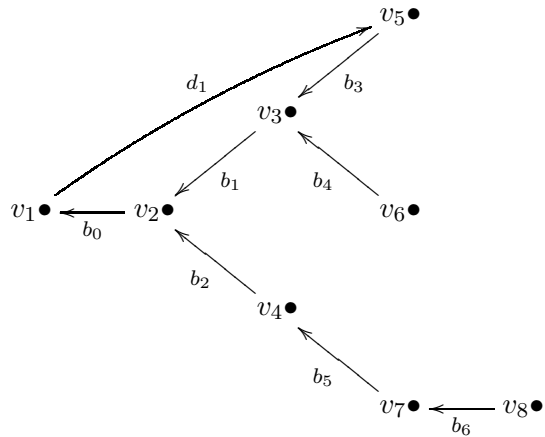
\xymatrix{
&&&v_5\bullet\ar[dl]^{b_3}&\\
&&v_3\bullet\ar[dl]^{b_1}&&\\
v_1\bullet\ar@/^/[rrruu]^{d_1}&v_2\bullet\ar[l]^{b_0}&&v_6\bullet\ar[ul]^{b_4}&\\
&&v_4\bullet\ar[ul]^{b_2}&&\\
&&&v_7\bullet\ar[ul]^{b_5}&v_8\bullet\ar[l]^{b_6}}
$$

We move on to the next level and conclude that  we need to add an
arrow from $v_2$ to $v_4$.  We do not add an arrow from $v_2$ to
$v_5$ because the arrow from $v_3$ to $v_2$  is already in a fully
recovered cycle. For level two vertices we need to add an arrow
from $v_3$ to $v_6$ and an arrow from $v_4$ to $v_8$.  Finally,
the recovered component is given by the quiver
$$
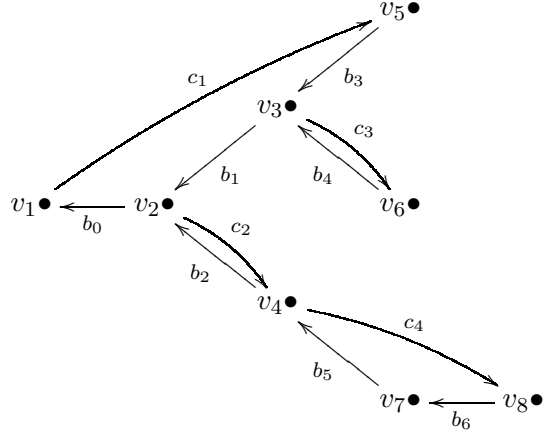
\xymatrix{
&&&v_5\bullet\ar[dl]^{b_3}&\\
&&v_3\bullet\ar@/^/[dr]^{c_3}\ar[dl]^{b_1}&&\\
v_1\bullet\ar@/^/[rrruu]^{c_1}&v_2\bullet\ar@/^/[dr]^{c_2}\ar[l]^{b_0}&&v_6\bullet\ar[ul]^{b_4}&\\
&&v_4\bullet\ar@/^/[rrd]^{c_4}\ar[ul]^{b_2}&&\\
&&&v_7\bullet\ar[ul]^{b_5}&v_8\bullet\ar[l]^{b_6}}
$$
}
\end{ex}

\begin{te}
Let $A_{\Gamma}$ be a graded basic Brauer tree algebra whose
grading is constructed by taking Green's walk around $\Gamma$.
From the quiver and relations of $A_0$ and the cyclic ordering of
the components of $A_0$ we can recover the quiver and relations of
$A_{\Gamma}$.
\end{te}
\noindent {\bf Proof.} We have seen in the previous proposition
that from the quiver and relations of $A_0$ we can recover each of
the components of the quiver of $A_{\Gamma}$ that we get when we
omit the exceptional cycle. In order to completely recover the
quiver of $A_{\Gamma}$, we are left to recover the exceptional
cycle. The roots of the components of $A_0$ are the vertices of
the exceptional cycle. From the cyclic ordering of the components
we get the cyclic ordering of the vertices of the exceptional
cycle. Thus, the exceptional cycle is recovered from the cyclic
ordering of the components of $A_0$. $\blacksquare$

\subsection{Quasi-hereditary structure on $ A_0$}
Let $Q_v$ be the quiver of an arbitrary connected component of
$A_0$. We have seen that $Q_v$ is a rooted tree. We can enumerate
the vertices of $Q_v$ in a natural way by the levels of the rooted
tree. We start with the root $v$, then we enumerate all  vertices
that are in level 1 of the rooted tree, for example, we enumerate
them from left to right in the planar embedding of the tree. Once
we have enumerated all vertices of an arbitrary level $r$, we move
on to level $r+1$ and repeat the same procedure until we enumerate
all vertices. Let $P_i$ be the projective cover of the simple
$A_0$-module $S_i$ corresponding to the vertex $v_i$. Then $P_i$
is spanned by paths of $Q_v$ ending at $v_i$. Since $Q_v$ is a
rooted tree, we conclude that the only simple modules that occur
as composition factors of $P_i$ are the simple modules whose
corresponding vertex has index greater than $i$. Also, $S_i$
occurs only once as a composition factor of $P_i$. Hence, we
obtain a quasi-hereditary structure on this component, by defining
a partial order as follows. Let $v_j$ be the vertex of $Q_v$
corresponding to the simple module $S_{j}$. Then we define $S_{j}<
S_{i}$, for $i\neq j$, if there is a path from $v_j$ to $v_i$,
where $S_{i}$ is the simple $A_0$-module corresponding to the
vertex $v_i$. The standard modules with respect to this order are
the projective indecomposable modules and  the costandard modules
are the simple modules. Therefore, $(A_0,\leq)$ is a
quasi-hereditary algebra as a product of quasi-hereditary
algebras.

The Cartan matrix of the path algebra of $Q_v$ is a lower
triangular matrix with diagonal elements equal to 1. Since $ A_0$
is the product of its components, we have that the following
standard result for quasi-hereditary algebras holds for $A_0$
(cf.\ [\ref{KenHer}]).

\begin{prop}
Let $\Gamma $ be a Brauer tree of type $(m,e)$ and let
$A_{\Gamma}$ be a graded basic Brauer tree algebra whose tree is
$\Gamma$ and whose grading is constructed by taking Green's walk
around $\Gamma$. If $A_0$ is the subalgebra of $A_{\Gamma}$
consisted of elements in degree $0$, then the Cartan matrix of
$A_0$ is a lower triangular matrix with diagonal elements equal to
$1$ and with determinant equal to $1$.
\end{prop}

Quasi-hereditary algebras have finite global dimension (cf.\
[\ref{RingelQHer}]), hence, $A_0$ has a finite global dimension.
We give an upper bound for the global dimension of a
quasi-hereditary algebra $A_0$. Let $Q_v$ be the quiver  of a
connected component of $A_0$ and let $B$ be its path algebra. Let
$l(Q_v)$ be the length of the rooted tree $Q_v$, ie.\ the index of
the last level of $Q_v$. The global dimension of $B$ is at most
$l(Q_v)$. This can be easily proved by looking at the projective
dimensions of the simple $B$-modules. One starts at the bottom of
the tree and works by induction on the distance of a vertex from
the bottom of the tree.

\begin{prop}
Let $A_{\Gamma}$, $A_0$ and $Q_v$ be as above. Then,
$${\rm gl.dim.}\, A_0\leq\max \{l(Q_v)\,|\, Q_v\,\, {\rm  a\,\, component\,\, of\,\, the \,\, quiver\,\, of}\,\, A_0\}.$$
\end{prop}

Note that the upper bound is achieved if the relations of $A_0$
are maximal possible in a sense that the product of every two
arrows is equal to zero. For example, this happens in the case of
a  Brauer line where the subalgebra $A_0$ is given by the quiver
$$
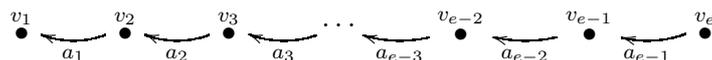
\xymatrix{\stackrel{v_1}{\bullet}&\stackrel{v_2}{\bullet}\ar@/^/[l]^{a_1}&\stackrel{v_3}{\bullet}\ar@/^/[l]^{a_2}
&\dots\ar@/^/[l]^{a_{3}}&\stackrel{v_{e-2}}{\bullet}\ar@/^/[l]^{a_{e-3}}&\stackrel{v_{e-1}}{\bullet}\ar@/^/[l]^{a_{e-2}}&\stackrel{v_e}{\bullet}\ar@/^/[l]^{a_{e-1}}}
$$
and the relations are $a_ia_{i-1}=0$, $i=2,3,\dots, e-1$. The
other extreme is the case when there are no relations, that is
when $A_0$ is hereditary. Then $A_0$ has global dimension $\leq
1$, on the other hand $l(Q_v)$ can be arbitrarily large.

\section{Graded Cartan matrix}
Let $A_{\Gamma}$ be a basic Brauer tree algebra of type $(m,e)$
given by the quiver $Q$ and relations $I$. We have seen that
$A_{\Gamma}$ is a graded algebra. Let $S_1, S_2,\dots , S_e$ be
the simple $A_{\Gamma}$--modules corresponding to the vertices of
the quiver $Q$. We assume that the simple modules are enumerated
by taking Green's walk around $\Gamma$. We define the graded
Cartan matrix $C$ of $A_{\Gamma}$ to be the ($e\times e$)-matrix
with entries from the ring $\mathbb{Z}[q,q^{-1}]$ given by
$$c_{ij}=C(S_i,S_j):=\sum_{l\in \mathbb{Z}}q^l {\rm dim}\, {\rm Homgr}_{A_{\Gamma}}(P_{S_i},P_{S_j}\langle l \rangle ),$$
where $P_{S_i}$ is the projective cover of $S_i$.

Note that the coefficient of $q^l$ is equal to the number of times
that $S_i$ appears in degree $l$ as a composition factor of
$P_{S_j}$.

\begin{prop}
Let $A_{\Gamma}$ be a graded basic Brauer tree algebra whose tree
is $\Gamma$ and with grading constructed by taking Green's walk
around $\Gamma$. Let $S_i$ and $S_j$ be simple modules
corresponding to vertices $v_i$ and $v_j$ of the quiver $Q$ of
$A_{\Gamma}$. Then
\begin{itemize}
\item[{\rm (i)}] if $S_i$ and $S_j$ do not belong to the same
cycle of $Q$,  then $c_{ij}=0$;

\item[{\rm (ii)}] if $S_i$ belongs to the exceptional cycle, we
have that
$$ c_{ii}=1+q^e+q^{2e}+\dots +
q^{me},$$   if $S_i$ does not belong to the exceptional cycle, we
have that
$$c_{ii}=1+q^{me}.$$

\item[{\rm (iii)}] if $i\neq j$ and $S_i$ and $S_j$ belong to the
same non-exceptional cycle, then
$$i>j \Rightarrow c_{ij}=1,$$
$$i<j \Rightarrow c_{ij}=q^{me};$$

\item[{\rm (iv)}] if $i\neq j$ and $S_i$ and $S_j$ belong to the
exceptional cycle, then
$$i>j \Rightarrow c_{ij}=q^{i-j}+q^{i-j+e}+\dots + q^{i-j+(m-1)e},$$
$$i<j \Rightarrow c_{ij}=q^{e-(j-i)}+q^{2e-(j-i)}+\dots + q^{me-(j-i)}.$$
\end{itemize}
\end{prop}

\noindent{\bf Proof.} Since the projective cover of $S_j$ is
spanned by the paths ending at $S_j$, we conclude that the
exponents of the non-zero terms of $c_{ij}$ are exactly the
degrees of the non-zero paths starting at $S_i$ and ending at
$S_j$. Case ${\rm (i)}$ is obvious, because $P_{S_j}$ does not
contain $S_i$ as a composition factor. In case ${\rm (ii)}$ the
degrees of the paths starting and ending at $S_i$ are $0, e, 2e,
\dots, me$ when $S_i$ belongs to the exceptional cycle, and are
$0,me$ otherwise. In case ${\rm (iii)}$, if $i>j$, the only
non-zero path from $S_i$ to $S_j$ has degree $0$. Similarly, if
$i<j$,  the only non-zero path from $S_i$ to $S_j$ has degree
$me$. In  case ${\rm (iv)}$  the same argument shows that, if
$i>j$, then the degrees of the paths from $S_i$ to $S_j$ are $i-j,
e+(i-j), \dots , (m-1)e+i-j$, and if $i<j$, they are $e-(j-i),
2e-(j-i), \dots , me-(j-i)$. $\blacksquare$

\section{Graded Cartan determinant}

\begin{prop}
Let $A_{\Gamma}$ be a graded basic Brauer tree algebra whose tree
 $\Gamma$ is of type $(m,e)$ and whose  grading is constructed by
taking Green's walk around $\Gamma$. If $C_{A_{\Gamma}}$ is the
graded Cartan matrix of $A_{\Gamma}$, then
$${\rm det}\, C_{A_{\Gamma}}=1+q^e+q^{2e}+\dots + q^{me^2}.$$
\end{prop}

\noindent{\bf Proof.} By [\ref{Rou}], Proposition 5.17, the
constant term of ${\rm det } \, C_A$ is equal to the determinant
of the Cartan matrix of $A_0$. We have seen that the determinant
of the Cartan matrix of $A_0$ is $1$. By Proposition 6.6 in
[\ref{Rou}], we also have that if $A$ and $B$ are two graded
Brauer tree algebras of the same type $(m,e)$, with gradings
constructed by taking Green's walk, then ${\rm det} \, C_A$ is
equal to $\pm q^l\,\, {\rm det}\, C_B$ for some integer $l$. Since
the constant term is equal to 1, we conclude that $l=0$, and that
${\rm det}\, C_A={\rm det}\, C_B$. Therefore, it is enough to
compute ${\rm det}\, C_B$ where $B$ is the graded basic Brauer
tree algebra whose tree is the Brauer line of type $(m,e)$ with
the exceptional vertex at one of the ends (see Example
\ref{linePR}(b)).

If $|i-j|>1$, then $c_{ij}=0$, because the corresponding vertices
belong to different cycles. Also, $c_{11}=1+q^e+q^{2e}+\dots
+q^{me} $, and if $i>1$, then $c_{ii}=1+q^{me}$. Other entries are
given by $c_{i,i+1}=q^{me}$ and $c_{i+1,i}=1$. We are left to
compute the following $e\times e$ determinant

$${\rm det}\, C_B=\left|
\begin{array}{cccccc}
  \alpha & \beta &  &  &  &  \\
  1 & \gamma & \beta &  &  &  \\
    &1 & \gamma & \beta &  & \\
    & &  & \ddots &  &  \\
   &  &  & 1 & \gamma & \beta \\
   &  & &  & 1 & \gamma \\
\end{array}\right|
$$
where $\alpha=1+q^e+q^{2e}+\dots +q^{me}$, $\beta=q^{me}$,
$\gamma=1+q^{me}$ and the omitted entries are all equal to zero.
If $d_l$ is the determinant of the  $l\times l$ block in the lower
right corner, then from $d_0=1$, $d_1=\gamma$ and the recursion
$$d_l=\gamma d_{l-1}-\beta d_{l-2},$$ it is easy to show that
$$d_l=1+q^{me}+\dots+q^{lme}.$$ Expanding the determinant along
the first column gives us the desired formula $${\rm det}\,
C_B=\alpha d_{e-1}-\beta d_{e-2}=1+q^e+q^{2e}+\dots+q^{me^2}.
$$ $\blacksquare$

\section{Brauer lines as trivial extension algebras}
Let $\Gamma$ be the  Brauer line with $e$ edges and multiplicity
of the exceptional vertex equal to 1. Let $A_{\Gamma}$ be a basic
Brauer tree algebra whose tree is $\Gamma$ and let us assume that
this algebra is graded by taking Green's walk around $\Gamma$. We
have seen in Example \ref{linePR} (b) that with respect to such
grading the graded quiver of $A_{\Gamma}$ is given by
$$
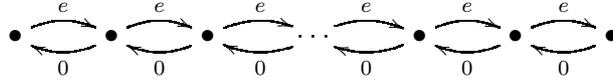
\xymatrix{\bullet\ar@/^/[r]^{e}&\bullet\ar@/^/[r]^{e}\ar@/^/[l]^0&\bullet\ar@/^/[l]^0\ar@/^/[r]^{e}&\dots\ar@/^/[l]^0\ar@/^/[r]^{e}&\bullet\ar@/^/[l]^0\ar@/^/[r]^{e}&\bullet\ar@/^/[l]^0\ar@/^/[r]^{e}&\bullet\ar@/^/[l]^0}$$
and we have that the only non-zero degree appearing in this
grading is $e$. Therefore, we can divide every degree by $e$ and
we will still have a graded algebra whose graded quiver is
$$
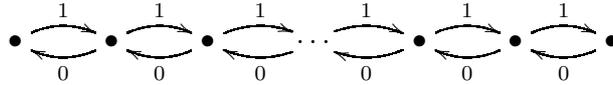
\xymatrix{\bullet\ar@/^/[r]^{1}&\bullet\ar@/^/[r]^{1}\ar@/^/[l]^0&\bullet\ar@/^/[l]^0
\ar@/^/[r]^{1}&\dots\ar@/^/[l]^0\ar@/^/[r]^{1}&\bullet\ar@/^/[l]^0\ar@/^/[r]^{1}&\bullet\ar@/^/[l]^0\ar@/^/[r]^{1}&\bullet\ar@/^/[l]^0}$$
with arrows only in degrees 0 and 1. We call this procedure of
dividing each degree by the same integer rescaling.

This algebra has an interesting connection with trivial extension
algebras. Let $B$ be a finite dimensional algebra over a field
$k$. The trivial extension algebra of $B$, denoted $T(B)$, is the
vector space $B\oplus B^*$ with multiplication defined by
$$(x,f)(y,g):=(xy,xg+fy)$$
where $x,y\in B$ and $f,g \in B^*$ and $B^*$ is the $B$--bimodule
${\rm Hom}_k(B,k)$. This algebra is always symmetric and the map
$B\rightarrow T(B)$, given by $b\mapsto (b,0)$, is an embedding of
algebras. The algebra $T(B)$ is naturally graded by putting $B$ in
degree 0 and $B^*$ in degree 1. This raises the question of
whether the graded Brauer tree algebra $A_{\Gamma}$ (with degrees
normalized by dividing by $e$) is the trivial extension algebra of
some algebra $B$? The obvious candidate would be its subalgebra
$A_0$. The quiver of $A_0$ is given by
$$
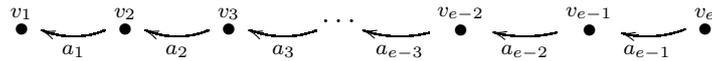
\xymatrix{\stackrel{v_1}{\bullet}&\stackrel{v_2}{\bullet}\ar@/^/[l]^{a_1}&\stackrel{v_3}{\bullet}\ar@/^/[l]^{a_2}
&\dots\ar@/^/[l]^{a_{3}}&\stackrel{v_{e-2}}{\bullet}\ar@/^/[l]^{a_{e-3}}&\stackrel{v_{e-1}}{\bullet}\ar@/^/[l]^{a_{e-2}}&\stackrel{v_e}{\bullet}\ar@/^/[l]^{a_{e-1}}}
$$
and the following proposition says that the trivial extension
algebra of $A_0$ is $A_{\Gamma}$.

\begin{prop}
Let $A_{\Gamma}$ and $A_0$ be as above. Then
$$T(A_0)=A_0\oplus A_0^*\cong  A_{\Gamma}.$$
\end{prop}
\noindent {\bf Proof.} Let $\{v_1^*,\dots, v_e^*, a_1^*,\dots,
a_{e-1}^*\}$ be the  basis of $A_0^*$ dual to the basis
$\{v_1,\dots, v_e, a_1,\dots, a_{e-1}\}$ of $A_0$ and let $b_i$,
$i=1,2,\dots, e-1$, be the arrow of the quiver of $A_{\Gamma}$
starting at the vertex $v_i$ and ending at the vertex $v_{i+1}$.
Each $b_i$ is in degree 1. It is now easily verified that the map
given by $a\mapsto (a,0)$ for $a \in A_0$ and $b_i\mapsto
(0,a_i^*)$, $i=1,2,\dots, e-1$, is an algebra isomorphism between
 $A_{\Gamma}$ and $T(A_0)$. $\blacksquare$

\section{Shifts of gradings}
Let $\Gamma $ be a Brauer tree of type $(m,e)$ and let
$A_{\Gamma}$ be a basic Brauer tree algebra whose tree is
$\Gamma$.  We have seen that we can grade this algebra by
computing the endomorphism ring of  the graded complex
$T^{\prime}=\oplus_{i=1}^{e} T_i^{\prime}$ which we constructed by
taking Green's walk around  $\Gamma$. In other words, we got a
structure of a graded algebra $A_{\Gamma}^{\prime}$ on
$A_{\Gamma}$. Recall that this was a non-negative grading, i.e.\ a
grading such that every homogeneous element is in a non-negative
degree. Let $\tilde{T}$ be the shifted graded complex
$\oplus_{i=1}^e T_i^{\prime}\langle n_i\rangle$, where $n_i\in
\mathbb{Z}$, $i=1,2,\dots, e$. The endomorphism ring of the graded
complex $\tilde{T}$ is the graded algebra $\tilde{A}_{\Gamma}$
which is graded Morita equivalent to $A_{\Gamma}^{\prime}$ (see
Definition \ref{grMorDef}). The question is if we can choose
non-zero integers $n_i$  in such way that the resulting grading is
positive, i.e.\ to get such a grading in which all homogeneous
elements from the radical of $A_{\Gamma}$ are in strictly positive
degrees. The answer to this question is positive, and moreover,
these integers $n_i$ can be chosen to be positive integers.

Let $S_i$ and $S_j$ be  vertices of the quiver of the algebra
$\tilde{A}_{\Gamma}$  which belong to the same cycle, and which
correspond to the summands $T_i^{\prime}\langle n_i\rangle$ and
$T_j^{\prime}\langle n_j\rangle$ of $\tilde{T}$. We need to
compute the degree of  the arrow $\alpha$  from $S_i$ to $S_j$.
Let $d$ be the degree of this arrow for the graded algebra
$A_{\Gamma}^{\prime}$.

\begin{prop}\label{shiftchange}
Let $\alpha$ be the arrow connecting vertices $S_i$ and $S_j$ of
the graded quiver of  the graded algebra $\tilde{A}_{\Gamma}$.
Then,
$${\rm deg }(\alpha)=d+n_i-n_j,$$
where $d$ is the degree of the same arrow of the quiver of the
graded algebra $A_{\Gamma}^{\prime}$, and $n_i$ and $n_j$ are the
shifts of $T_i^{\prime}$ and $T_j^{\prime}$ respectively.
\end{prop}
\noindent {\bf Proof.} Since $T_i^{\prime}$ are complexes of
uniserial modules,  the top of the leftmost non-zero term of
$T_i^{\prime}\langle n_i\rangle$ is in  degree $d_1-n_i$ after the
shift, where $d_1$ is the degree of the top of the leftmost
non-zero term of $T_i^{\prime}$, ie.\ the degree before the shift.
Also, the top of the leftmost non-zero term  of
$T_j^{\prime}\langle n_j\rangle$ is in degree $d_2-n_j$, where
$d_2$ is its degree before the shift. The degree of $\alpha$ after
the shift is $(d_2-n_j)-(d_1-n_i)=d+n_i-n_j$. $\blacksquare$

Note that when we compare  graded quivers of $A_{\Gamma}^{\prime}$
and $\tilde{A}_{\Gamma}$, the difference is that we added
$n_i-n_j$ to the degree of the arrow from $S_i$ to $S_j$. If this
arrow is in degree 0, then its degree after the shift is
$n_i-n_j$, where $i>j$. The source and the target of such an arrow
belong to two consecutive levels of a rooted tree. The arrows that
are not part of the exceptional cycle and whose degree was
non-zero, are now in degree $me+n_i-n_j$ where $i<j$. Also, the
degree of an arrow between two vertices $S_i$ and $S_j$ of the
exceptional cycle is now $i-j+n_i-n_j$, if $i>j$, and if $i<j$, it
is $e-(j-i)+n_i-n_j$.

We want to find integers $n_i$, $i=1,2,\dots, e$, such that all
these degrees are positive.

\begin{prop}
Let $\Gamma$ be a Brauer tree of type $(m,e)$ and let
$\tilde{T}:=\oplus_{i=1}^e T_i^{\prime}\langle n_i\rangle $ be the
shifted tilting complex constructed by taking Green's walk around
$\Gamma$. Let $\tilde{A}_{\Gamma}:={\rm Endgr}_{K^b(P_{A_S})}
(\tilde{T})^{op}$. There are positive integers $n_i$,
$i=1,2,\dots, e$, such that the graded algebra
$\tilde{A}_{\Gamma}$ is non-negatively graded with ${\rm
deg}(a)>0$ for all homogeneous elements $a\in {\rm rad}\,
\tilde{A}_{\Gamma}$.
\end{prop}
\noindent{\bf Proof.} Let $Q_v$ be an arbitrary component of the
quiver of $A_0$, where $A_0$ is as before the subalgebra of
$A^{\prime}_{\Gamma}$ of degree 0 elements, and let $S_i$ be a
vertex of $Q_v$. If we choose $n_i$ to be $1+l_i$, where $l_i$ is
the level of the rooted tree $Q_v$ to whom $S_i$ belongs, we see
that all arrows of the graded quiver $Q_v$ are in degree 1 after
the shift. Also, the arrows of $Q$, the quiver of $A_{\Gamma}$,
that connect two vertices of $Q_v$ and which were not in degree 0
are still in positive degrees after the shift because $n_j-n_i<e$
(the number of levels in each component of the quiver of $A_0$ is
less than $e$) and consequently $me+n_i-n_j>0$ . The arrows of the
exceptional cycle are in the same degrees as they used to be
because we set $n_i:=1$ for each root $S_i$. Then for every
homogeneous element $a\in {\rm rad}\, \tilde{A}_{\Gamma}$ we have
that ${\rm deg}(a)>0$. $\blacksquare$

We note here that, in general, there are many choices for the
integers $n_i$, $i=1,2,\dots,e$.

\begin{ex}
{\rm  Let $\Gamma$ be the tree from  Example \ref{NulaPr}. With
the notation from the previous proposition the graded quiver of
the basic Brauer tree algebra $\tilde{A_{\Gamma}}={\rm
Endgr}_{K^b(P_{A_S})}(\oplus_{i=1}^{11}T_i^{\prime}\langle
n_i\rangle)^{op}$ is given by
$$
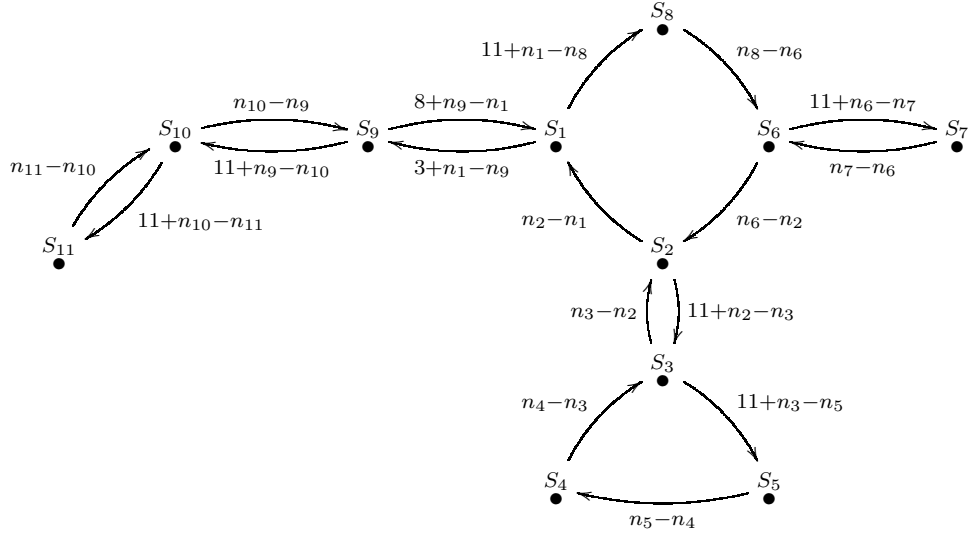
\xymatrix{
&&&&&&\stackrel{S_8}{\bullet}\ar@/^/[dr]^{n_8-n_6}&&\\
&\stackrel{S_{10}}{\bullet}\ar@/^/[dl]^{11+n_{10}-n_{11}}\ar@/^/[rr]^{n_{10}-n_9}&&\stackrel{S_9}{\bullet}\ar@/^/[ll]^{11+n_9-n_{10}}\ar@/^/[rr]^{8+n_9-n_1}&&\stackrel{S_1}{\bullet}\ar@/^/[ll]^{3+n_1-n_9}\ar@/^/[ur]^{11+n_1-n_8}&&\stackrel{S_6}{\bullet}\ar@/^/[dl]^{n_6-n_2}\ar@/^/[rr]^{11+n_6-n_7}&&\stackrel{S_7}{\bullet}\ar@/^/[ll]^{n_7-n_6}\\
\stackrel{S_{11}}{\bullet}\ar@/^/[ur]^{n_{11}-n_{10}}&&&&&&\stackrel{S_2}{\bullet}\ar@/^/[d]^{11+n_2-n_3}\ar@/^/[ul]^{n_2-n_1}&&&\\
&&&&&&\stackrel{S_3}{\bullet}\ar@/^/[u]^{n_3-n_2}\ar@/^/[dr]^{11+n_3-n_5}&&&\\
&&&&&\stackrel{S_4}{\bullet}\ar@/^/[ur]^{n_4-n_3}&&\stackrel{S_5}{\bullet}\ar@/^/[ll]^{n_5-n_4}&&
}$$ If we set $n_1=n_9=1$, $n_2=n_{10}=2$, $n_3=n_6=n_{11}=3$,
$n_4=n_7=n_8=4$ and $n_5=5$, then all arrows are in positive
degrees. }
\end{ex}

\noindent Note that the change of shifts on the summands
$T_i^{\prime}$ of the tilting complex $T^{\prime}$ is the same as
the change of shifts on the projective indecomposable modules of
$A_{\Gamma}^{\prime}\cong \displaystyle{\rm
Endgr}_{K^b(P_{A_S})}(\oplus_{i=1}^e T_i^{\prime})^{op}$. Let
$\tilde{A}_{\Gamma}\cong  {\rm
Endgr}_{K^b(P_{A_S})}(\oplus_{i=1}^e T_i^{\prime}\langle
n_i\rangle)^{op}$. When we change the shifts, in general, we get a
different grading on $A_{\Gamma}$ and the resulting graded algebra
$\tilde{A}_{\Gamma}$ is not isomorphic to $A_{\Gamma}^{\prime}$ as
a graded algebra. But these two graded algebras are graded Morita
equivalent, ie.\ there is an equivalence
$A_{\Gamma}^{\prime}$--${\rm grmod}\cong
\tilde{A}_{\Gamma}$--${\rm grmod}$ as we shall see in Section 11.

\section{Change of the exceptional vertex}
Let $A_{\Gamma}$ be a basic Brauer tree algebra whose tree
$\Gamma$ has $e$ edges and  the multiplicity of the exceptional
vertex equal to 1. If we change the exceptional vertex, the
algebra $A_{\Gamma}$ does not change. But when constructing the
tilting complex that tilts from $A_S$ to $A_{\Gamma}$ by taking
Green's walk around $\Gamma$ it is obvious that we start from a
different vertex, and in general, the resulting tilting complex is
different. Therefore, we get different gradings on  $A_{\Gamma}$.

\begin{ex}\label{promcvora}
{\rm  Let $\Gamma$ be the following Brauer tree with multiplicity
of the exceptional vertex equal to 1 and let $A_{\Gamma}$ be the
corresponding  Brauer tree algebra.

$$
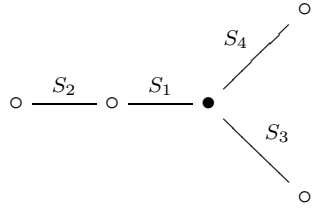
\xymatrix{&&& \circ\ar@{-}[dl]_{S_4}\\
\circ\ar@{-}[r]^{S_2}&\circ\ar@{-}[r]^{S_1}&\bullet\ar@{-}[dr]^{S_3}&\\
&&&\circ\\}$$

If $T=\oplus_{i=1}^4T_i$ is the tilting complex constructed by
taking Green's walk around $\Gamma$, then the resulting graded
quiver of  $A_{\Gamma}$ is given by

$$
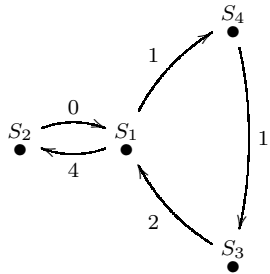
\xymatrix{&&\stackrel{S_4}{\bullet}\ar@/^/[dd]^1\\
\stackrel{S_2}{\bullet}\ar@/^/[r]^0&\stackrel{S_1}{\bullet}\ar@/^/[ur]^1\ar@/^/[l]^4&\\
&&\stackrel{S_3}{\bullet}\ar@/^/[ul]^2}
$$

If we change the exceptional vertex, say we have Brauer tree
$\Delta$

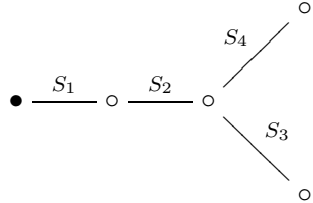
$$\xymatrix{&&& \circ\ar@{-}[dl]_{S_4}\\
\bullet\ar@{-}[r]^{S_1}&\circ\ar@{-}[r]^{S_2}&\circ\ar@{-}[dr]^{S_3}&\\
&&&\circ\\}$$

\noindent then the basic Brauer tree algebra $A_{\Delta}$ whose
tree is  $\Delta$, is the same as $A_\Gamma$. Thus, changing the
exceptional vertex of $\Gamma$ does not change $A_{\Gamma}$. The
tilting complex $D$ constructed by taking Green's walk around
$\Delta$ is different from $T$. Therefore, we get  a new grading
on $A_{\Delta}=A_{\Gamma}\cong {\rm End}_{K^b(P_{A_S})}(D)^{op}$,
and the resulting graded quiver of the graded algebra
$A_{\Delta}^{\prime}$ is given by

$$
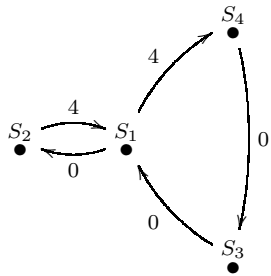
\xymatrix{&&\stackrel{S_4}{\bullet}\ar@/^/[dd]^0\\
\stackrel{S_2}{\bullet}\ar@/^/[r]^4&\stackrel{S_1}{\bullet}\ar@/^/[ur]^4\ar@/^/[l]^0&\\
&&\stackrel{S_3}{\bullet}\ar@/^/[ul]^0}
$$

\noindent If $\tilde{T}:=\oplus_{i=1}^4T_i^{\prime}\langle
n_i\rangle$ is a graded complex given by shifting summands of
$T^{\prime}$, then from Proposition \ref{shiftchange} we get
another grading on $A_{\Gamma}\cong {\rm
End}_{K^b(P_{A_S})}(T)^{op}$ and the resulting graded quiver of
the graded algebra $\tilde{A}_{\Gamma}$ is given by

$$
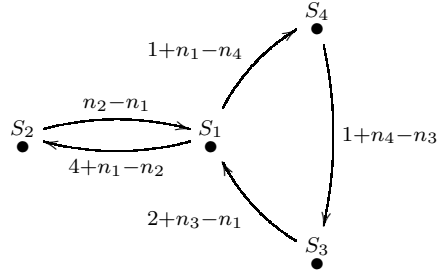
\xymatrix{&&&\stackrel{S_4}{\bullet}\ar@/^/[dd]^{1+n_4-n_3}\\
\stackrel{S_2}{\bullet}\ar@/^/[rr]^{n_2-n_1}&&\stackrel{S_1}{\bullet}\ar@/^/[ur]^{1+n_1-n_4}\ar@/^/[ll]^{4+n_1-n_2}&\\
&&&\stackrel{S_3}{\bullet}\ar@/^/[ul]^{2+n_3-n_1}}
$$

If we set $n_1=3$, $n_2=7$, $n_3=1$ and $n_4=0$, we see that the
resulting grading on $A_{\Gamma}$ is the same as the grading that
we got by taking Green's walk around $\Delta$. }
\end{ex}

In the previous example we had two different gradings on
$A_{\Gamma}$, but we were able, by changing shifts of the summands
of the graded tilting complex $T^{\prime}$, to move from one
grading to another via graded Morita equivalence (Definition
\ref{grMorDef}). We will prove in the next section that this holds
for all Brauer tree algebras, regardless of the multiplicity of
the exceptional vertex.

\section{Classification of gradings}

In this section we classify, up to graded Morita equivalence and
rescaling, all gradings on an arbitrary Brauer tree algebra with
$n$ edges and the multiplicity of the exceptional vertex equal to
$m$.

For a finite dimensional $k$-algebra $A$, there  is a
correspondence between gradings on $A$ and homomorphisms of
algebraic groups from $\textbf{G}_m$ to ${\rm Aut}(A)$, where
$\textbf{G}_m$ is the multiplicative group $k^*$ of a field $k$.
For each grading $A=\oplus_{i\in \mathbb{Z}}A_i$ there is a
homomorphism of algebraic groups $\pi \, : \, \textbf{G}_m
\rightarrow {\rm Aut}(A)$ where element $x\in k^*$ acts on $A_i$
by multiplication by $x^i$ (see [\ref{Rou}], Section 5). If $A$ is
graded and $\pi$ is the corresponding homomorphism, we will write
$(A,\pi)$ to denote that $A$ is graded with grading $\pi$.

\begin{de}\label{grMorDef} Let $(A,\pi)$ and $(A,\pi^{\prime})$ be two gradings
on a $k$-algebra $A$,  and let $P_1, P_2,\dots, P_r$ be the
isomorphism classes of the projective indecomposable $A$-modules.
We say that $(A,\pi)$ and $(A,\pi^{\prime})$ are graded Morita
equivalent if there exist integers $d_1,d_2,\dots,d_r$  such that
the graded algebras $(A,\pi^{\prime})$ and ${\rm
Endgr}_{(A,\pi)}(\oplus_{i=1}^rn_iP_i\langle d_i\rangle)^{\rm op}$
are isomorphic.
\end{de}
Note that two graded algebras are graded Morita equivalent if and
only if their categories of graded modules are equivalent.

The following proposition tells us how to classify all gradings on
$A$ up to graded Morita equivalence.

\begin{prop}[[\ref{Rou}{]}, Corollary 5.8]\label{cokarak}
Two graded algebras $(A,\pi)$ and $(A,\pi^{\prime})$ are graded
Morita equivalent if and only if the corresponding cocharacters
$\,\pi \, : \, \textbf{G}_m \rightarrow {\rm Out}(A)$  and
$\,\pi^{\prime} \, : \, \textbf{G}_m \rightarrow {\rm Out}(A)$ are
conjugate.
\end{prop}
From this proposition we see that in order to classify gradings on
$A$ up to graded Morita equivalence, we need to compute maximal
tori in  ${\rm Out}(A)$. Let ${\rm Out}^K (A)$ be  the subgroup of
${\rm Out}\, (A)$ of those automorphisms fixing the isomorphism
classes of simple $A$-modules. Since ${\rm Out}^K (A)$ contains
${\rm Out}^0(A)$, the connected component of ${\rm Out}(A)$ that
contains the identity element, we have that maximal tori in ${\rm
Out}(A)$ are actually contained in ${\rm Out}^K (A)$. \vspace{3mm}

We start by computing ${\rm Out}^K (A_{\Gamma})$ for a Brauer tree
algebra $A_{\Gamma}$ of type $(m,n)$, where $m>1$. We remark here
that for the case $m=1$, ${\rm Out}^K (A_{\Gamma})$ has been
computed in Lemma 4.3 in [\ref{RouZim}]. Although the proof of
this Lemma in [\ref{RouZim}] seems to be incomplete, the result is
correct if one assumes that the ground field $k$ is an
algebraically closed field. The same result, namely that ${\rm
Out}^K (A_{\Gamma})\cong k^*$ when $m=1$, follows directly from
our computation below for $m>1$, if we disregard the loop $\rho$
which does not appear in the case $m=1.$

\vspace{3mm} Let $A_{\Gamma}$ be a basic Brauer tree algebra whose
tree $\Gamma$ is of type $(m,n)$ where $m>1$. Since ${\rm Out}^K
(A_{\Gamma})$ is invariant under derived equivalence (cf.\
[\ref{Link}], Section 4), in order to compute this group for any
Brauer tree algebra $A_{\Gamma}$ whose tree is of type $(m,n)$, it
is sufficient to compute this group for the Brauer line of the
same type.

Let $A_{\Gamma}$ be a basic Brauer tree algebra of type $(m,n)$
whose tree is the Brauer line with $n$ edges and with the
exceptional vertex at the end of the line. The quiver of
$A_{\Gamma}$ is given by
$$
\xymatrix{\stackrel{e_1}{\bullet}\ar@(ul,dl)_{\rho}\ar@/^/[r]^{a_1}&\stackrel{e_2}{\bullet}\ar@/^/[r]^{a_2}\ar@/^/[l]^{b_1}&\stackrel{e_3}{\bullet}\ar@/^/[l]^{b_2}\ar@/^/[r]^{a_3}&\dots\ar@/^/[l]^{b_3}\ar@/^/[r]^{a_{n-3}}&\stackrel{e_{n-2}}{\bullet}\ar@/^/[l]^{b_{n-3}}\ar@/^/[r]^{a_{n-2}}&\stackrel{e_{n-1}}{\bullet}\ar@/^/[l]^{b_{n-2}}\ar@/^/[r]^{a_{n-1}}&\stackrel{e_n}{\bullet}\ar@/^/[l]^{b_{n-1}}}
$$
and relations are given by $a_ia_{i+1}=b_{i+1}b_i=0$
$(i=1,2,\dots, n-2)$, $\rho^m=a_1b_1$ and $a_ib_i=b_{i-1}a_{i-1}$
$(i=2,3,\dots, n-1)$.

In order to simplify calculation let us set $t_i:=a_i+b_i$,
$i=1,2,\dots,n-1.$ Then, $a_i=e_it_i$ and $b_i=e_{i+1}t_i$, for
$i=1,2,\dots, n-1$. Then the  relations become $t_i^3=0$
$(i=1,2,\dots,n-1)$, $e_1t_1^2=\rho^m$, $e_it_i^2=t_{i-1}^2e_i$
$(i=2,3,\dots,n-1)$. Also, we have that $e_it_i=t_ie_{i+1}$ and
$e_{i+1}t_i=t_ie_i$,  $i=1,2,\dots, n-1$. Let us write
$t_0^2:=\rho^m$.

Let $\varphi$ be an arbitrary automorphism of $A_{\Gamma}$ that
fixes isomorphism classes of simple $A_{\Gamma}$-modules. We will
compose this automorphism with suitably chosen inner automorphisms
in order to get an automorphism $\phi$, which represents the same
element as $\varphi$ in the group of outer automorphisms, but
which is suitable for our computation.

If $\{e_1, e_2,\dots, e_n\}$ is a complete set of orthogonal
primitive idempotents, then
$\{\varphi(e_1),\varphi(e_2),\dots,\varphi(e_n)\}$ is also a
complete set of orthogonal primitive idempotents. From classical
ring theory (cf.\ [\ref{jac}], Theorem 3.10.2) we know that these
two sets are conjugate. Hence, when we compose $\varphi$ with a
suitably chosen inner automorphism we get that
$x^{-1}\varphi(e_i)x=e_ {\pi(i)}$, for all $i$, where $\pi$ is
some permutation. Since $\varphi$ fixes isomorphism classes of
simple modules, we can assume that, for all $i$,
$\varphi(e_i)=e_i.$

Since $\varphi({\rm rad}\, A_{\Gamma})\subset {\rm rad}\,
A_{\Gamma}$ we have that
$$\varphi(t_i)=\sum_{j=1}^{n-1}A_{ij}e_jt_j+\sum_{j=1}^{n-1}B_{ij}e_{j+1}t_j+\sum_{j=0}^{n-1}C_{ij}e_{j+1}t_j^2+\sum_{j=1}^{m-1}D_{ij}\rho^j,$$
where  $A_{ij}$, $B_{ij}$, $C_{ij}$ and $D_{ij}$ are scalars.  If
$i>1$, then from $e_1t_i=0$ we get that $D_{ij}=0$ for
$j=1,2,\dots,m-1$. If $i=1$, then from $e_1t_1=t_1e_2$ we get that
$D_{1j}=0$, $j=1,2,\dots,m-1$. In both cases we have that
$D_{ij}=0$ for all $i$ and all $j$.

If $l\notin \{i,i+1\}$, then $\varphi(e_lt_i)=0$ and
$\varphi(t_ie_l)=0$. From the first equality we get that
$A_{il}=B_{il-1}=0$ and from the second equality we get that
$B_{il}=A_{il-1}=0$. Subsequently, we have that $C_{il}=0$ for
$l\notin \{i,i-1\}$. Therefore, we must have that
$$\varphi(t_i)=A_{ii}e_it_i+B_{ii}e_{i+1}t_i+C_{ii-1}e_{i}t_{i-1}^2+C_{ii}e_{i+1}t_i^2.$$

From $e_it_i=t_ie_{i+1}$ it follows that $\varphi(e_i)
\varphi(t_i)=\varphi(t_i)\varphi(e_{i+1})$, and we have
$$C_{ii-1}=C_{ii}=0.$$

Also, $A_{ii}\neq 0$ and $B_{ii}\neq 0$. If one of them would be
zero, then $\varphi(t_i^2)=0$, which is in contradiction with
$\varphi$ being an automorphism.

Assume now that
$$\varphi(\rho)=\sum_{j=1}^{m-1}D_{j}\rho^j+\sum_{j=1}^{n-1}E_{j}e_jt_j+\sum_{j=1}^{n-1}F_{j}e_{j+1}t_j+\sum_{j=0}^{n-1}G_{j}e_{j+1}t_j^2,$$
for some scalars $D_j$, $E_j$, $F_j$ and $G_j$.

For $l>2$, from $e_l\rho=0$ we have that
$\varphi(e_l)\varphi(\rho)=0$ and it follows that $E_l=F_{l-1}=0$,
for $l=2,3,\dots,m-1$. Similarly, from $\rho e_l=0$ we have that
$F_l=E_{l-1}=0$, for $l=2,3,\dots,m-1$ and consequently, we have
that $G_j=0$ for all $j$. Therefore, we get
$$\varphi(\rho)=\sum_{j=1}^m D_j \rho^j.$$
Note that $D_1\neq 0$, because if $D_1=0$, then
$\varphi(\rho^m)=0$, which again contradicts  the fact that
$\varphi$ is an automorphism.

Gathering all the information we got, we conclude that, up to
inner automorphism, an arbitrary automorphism that fixes
isomorphism classes of simple $A_{\Gamma}$-modules has the
following action on a set of generators $\{e_i,\, t_l, \, \rho\,
|\, i=1,2,\dots,n ; \, l=1,2,\dots,n-1\}$
$$\varphi(e_i)=e_i,\quad \varphi(t_l)=A_{ll}e_lt_l+B_{ll}e_{l+1}t_l, \quad \varphi(\rho)=\sum_{j=1}^m D_j \rho^j.$$

\noindent To compute ${\rm Out}^K ( A_{\Gamma})$, for each
automorphism $\varphi$  we need to find an automorphism which is
in the same class as $\varphi$ in ${\rm Out}^K (A_{\Gamma})$, and
which acts by scalar multiplication on as many of the above
generators as possible. In order to do that we will compose
$\varphi$ with suitably chosen inner automorphisms.

First of all, let us see how an arbitrary inner automorphism acts
on our set of generators.

Let $y$ be an arbitrary invertible element in $A_{\Gamma}$. Then,
$$y=\sum_{j=1}^nl_je_j+\sum_{j=1}^{n-1}s_je_jt_j+\sum_{j=1}^{n-1}r_je_{j+1}t_j+\sum_{j=1}^{n-1}p_je_{j+1}t_j^2+\sum_{j=1}^mq_j\rho^j,$$
where $l_i\neq 0$, $i=1,\dots,n$. From $yy^{-1}=1$, we easily
compute scalars $s_j^{\prime}$, $r_j^{\prime}$, $p_j^{\prime}$,
$q_j^{\prime}$, in
$$y^{-1}=\sum_{j=1}^nl_j^{-1}e_j+\sum_{j=1}^{n-1}s_j^{\prime}e_jt_j+\sum_{j=1}^{n-1}r_j^{\prime}e_{j+1}t_j+\sum_{j=1}^{n-1}p_j^{\prime}e_{j+1}t_j^2+\sum_{j=1}^mq_j^{\prime}\rho^j.$$

\noindent The inner automorphism given by $y$  has the following
action on $\rho$

$$f_y(\rho):=y\rho y^{-1}=(l_1\rho+\sum_{j=1}^m q_j\rho^{j+1})y^{-1}=$$
$$=\rho+\sum_{j=1}^ml_1q_j^{\prime}\rho^{j+1}+\sum_{j=1}^ml_1^{-1}q_j\rho^{j+1}+\sum_{j=1}^m q_j \rho^{j+1} \sum_{j=1}^m q_j^{\prime}\rho^{j}=$$
$$=\rho+\rho(\sum_{j=1}^ml_1q_j^{\prime}\rho^{j}+\sum_{j=1}^ml_1^{-1}q_j\rho^{j}+\sum_{j=1}^m q_j \rho^{j} \sum_{j=1}^m q_j^{\prime}\rho^{j}+s_1r_1^{\prime}\rho^m)=\rho.$$

\noindent Therefore, an arbitrary inner automorphism fixes $\rho$.
From now on, we will use inner automorphisms given by elements of
the form
$$x=\sum_{j=1}^nl_je_j.$$
They will be enough to get a class representative of $\varphi$ in
${\rm Out}^K (A_{\Gamma})$ that is easy to work with.

If $f_x$ is the inner automorphism given by an invertible element
$x$, then we easily get that
$$f_x(t_i)=l_il^{-1}_{i+1}e_it_i+l_i^{-1}l_{i+1}e_{i+1}t_i.$$

\noindent When we compose $f_x$ and $\varphi$ we get that
$$f_x\circ \varphi(\rho)=\sum^{m}_{j=1}D_j\rho^j,$$
$$ f_x\circ \varphi(t_i)=A_{ii}l_il_{i+1}^{-1}e_it_i+B_{ii}l_i^{-1}l_{i+1}e_{i+1}t_i,$$
\noindent for all $i$, and
$$f_x\circ \varphi(e_i)=e_i.$$

\noindent We want to choose $l_i$'s so that we get
$$A_i:=A_{ii}l_il_{i+1}^{-1}=B_{ii}l_i^{-1}l_{i+1}$$
for $i=1,2,\dots,n-1.$ To do this we need to choose $l_i$'s in
such way that the following equality holds for all $i$
$$A_{ii}B_{ii}^{-1}=l_i^{-2}l_{i+1}^2.$$ We can choose $l_1=1$ and
then inductively, assuming that we have chosen
$l_1,l_2,\dots,l_i$, because we are working over an algebraically
closed field, we get $l_{i+1}$ from
$A_{ii}B_{ii}^{-1}l_i^{2}=l_{i+1}^2.$

If we choose such $x$, then the map $\varphi_1:=f_{x}\circ
\varphi$ has the following action on our generating set
$$\varphi_1(e_i)=e_i,\quad \varphi_1(\rho)=\sum_{j=1}^mD_j\rho^j,\quad \varphi_1(t_i)=A_{i}t_i.$$

\noindent From the relations $e_it_{i-1}^2=t_{i}^2e_{i}$, for
$i=2,3,\dots,n-1$, we get that $A_1^2=A_2^2=\dots=A_{n-1}^2.$ We
can assume that $A_1=A_2=\dots=A_{n-1}$, because if not, then by
multiplying $\varphi_1$ by an inner automorphism given by
$x_1=\sum_{i=1}^nr_ie_i$, where we set $r_1=1$ and then
inductively, $r_{i+1}=-r_i$ if $A_{i+1}=-A_i$, and $r_{i+1}=r_i$
if $A_{i+1}=A_i$, we get a new automorphism $\varphi_2$ such that
$\varphi_2(t_i)=A_1t_i$, for all $i$. Also from the relation
$\rho^m=e_1t_1^2$ we get that $A_1^2=D_1^m.$ This means that for a
fixed $D_1$ we have two choices for $A_1$, since there are two
square roots of $D_1^m$. These two values of $A_1$ will give us
two different automorphisms, but as before, we can assume that
after multiplying by an appropriate inner automorphism these two
automorphisms represent the same automorphism in ${\rm
Out}^K(A_{\Gamma})$.

We started with an arbitrary automorphism $\varphi$ that fixes the
isomorphism classes of simple $A_{\Gamma}$-modules and we showed
that in the group ${\rm Out}^K (A_{\Gamma})$ it represents the
same class as the element $\phi$ whose action is given  by
$$ \phi(e_i)=e_i, \quad \phi(\rho)=\sum^m_{j=1}D_j\rho^j, \quad \phi(t_s)=A_1 t_s,$$
where  $A_1^2=D_1^m$, $i=1,\dots,n$, and $s=1,\dots,n-1.$

Therefore, every element in ${\rm Out}^K (A_{\Gamma})$ is uniquely
determined by its action on $\rho$, that is, it is uniquely
determined by an $m$-tuple $(D_1,D_2,\dots,D_m)$ where $D_1\neq
0$. The map $\theta$ that assigns to each $m$-tuple
$D=(D_1,D_2,\dots, D_m)$ an isomorphism $\phi_D$, where
$\phi_D(\rho)=\sum_{j=1}^m D_j \rho^j$, is an isomorphism of
groups. But what is the group structure on the set $k^*\times
\underbrace{k\times k \times\dots\times k}_{m-1}$, where $k^*$
denotes non-zero elements of $k$? If $\alpha=(\alpha_1, \alpha_2,
\dots, \alpha_m)$ and $\beta=(\beta_1,\beta_2,\dots,\beta_m)$ are
two $m$-tuples and $\phi_{\alpha}$ and $\phi_{\beta}$ are two
corresponding automorphisms, then $\phi_{\beta}\circ
\phi_{\alpha}=\phi_{\beta*\alpha}$  gives us the definition of the
group operation $*$ on $k^*\times \underbrace{k\times k
\times\dots\times k}_{m-1}$. Computing $\phi_{\beta}\circ
\phi_{\alpha}$ gives us that

\begin{equation} \label{eqGrupa}
\beta * \alpha:=\left(\sum_{i=1}^l\alpha_i(\sum_{\tiny \begin{array}{c}k_1+\dots+k_i=l\\
k_1,\dots,k_i>0\end{array}}\beta_{k_1}\beta_{k_2}\cdots\beta_{k_i})\right)_{l=1}^m
\end{equation}

\noindent Here are first few coordinates explicitly $$\beta *
\alpha=(\alpha_1\beta_1,\, \alpha_1\beta_2+\alpha_2\beta_1^2,
\alpha_1\beta_3+2\alpha_2\beta_1\beta_2+\alpha_3\beta_1^3,\dots).$$

\begin{de} We define $H_m$ to be the group $(k^*\times \underbrace{k\times k \times\dots\times
k}_{m-1}\, , \, *\, )$, where the multiplication $*$ is given by
the above equation \eqref{eqGrupa}.
\end{de}

\noindent The identity element of $H_m$ is $(1,0,\dots,0)$, and
this element corresponds to the class of inner automorphisms. The
inverse element of an arbitrary $m$-tuple is easily computed
inductively from the definition of $*$. Associativity is verified
after elementary, but tedious computation.

\begin{lm}
The group $H_m$ is isomorphic to the group of automorphisms of the
polynomial  algebra $k[x]/(x^{m+1}).$
\end{lm}
\noindent {\bf Proof.} An arbitrary automorphism $f$ from ${\rm
Aut} (k[x]/(x^{m+1}))$ is given by its action on $x$. Since it has
to be surjective, and $f(x)^{m+1}=0$, we have that $x$ has to be
mapped to a polynomial $d_1x+d_2x^2+\dots+d_mx^m$, where $d_1\neq
0$. Therefore, every automorphism of ${\rm Aut}(k[x]/(x^{m+1}))$
is given by a unique $m$-tuple $(d_1,d_2,\dots,d_m)$ where
$d_1\neq 0$. The structure of a group on the set of all such
$m$-tuples is the same as for the group $H_m$. $\blacksquare$

Since the group ${\rm Out}^K(A)$ is invariant under derived
equivalence and Brauer tree algebras of the same type are derived
equivalent, we have the following theorem.

\begin{te}
Let $\Gamma$ be a Brauer tree of type $(m,n)$ and let $A_{\Gamma}$
be a basic Brauer tree algebra whose tree is $\Gamma$. Then
$${\rm Out}^K (A_{\Gamma})\cong {\rm Aut}(k[x]/(x^{m+1}))\cong H_m.$$
\end{te}

We see that the group of outer automorphisms that fix isomorphism
classes of simple modules of an arbitrary Brauer tree algebra
depends only on the multiplicity $m$ of the exceptional vertex and
does not depend on the number of edges. If we take an arbitrary
Brauer tree $\Gamma$, and if $\Gamma^{\prime}$ is any connected
Brauer subtree of $\Gamma$ that contains the exceptional vertex,
then the corresponding Brauer tree algebras $A_{\Gamma}$ and
$A_{\Gamma^{\prime}}$ have isomorphic groups of outer
automorphisms that fix isomorphism classes of simple modules. If
we take the subtree ${\Gamma}^{\prime}$ to be the exceptional
vertex with one edge adjacent to it, we get
$A_{\Gamma^{\prime}}=k[x]/(x^{m+1}).$

\begin{co}
Let $\Gamma$ be a Brauer tree of type $(m,n)$ and let
$\Gamma^{\prime}$ be an arbitrary connected Brauer subtree of
$\Gamma$ that contains the exceptional vertex. If $A_{\Gamma}$ and
$A_{\Gamma^{\prime}}$ are two basic Brauer tree algebras whose
trees are $\Gamma$ and $\Gamma^{\prime}$ respectively, then
$${\rm Out}^K (A_{\Gamma})\cong {\rm Out}^K(A_{\Gamma^{\prime}})\cong {\rm Aut}(k[x]/(x^{m+1})).$$
\end{co}

Let $L$ be the subgroup of $H_m$ consisting of the elements of the
form $(1,\alpha_2,\dots, \alpha_m)$ and let $K$ be the subgroup of
$H_m$ consisting of the elements of the form
$(\alpha_1,0,\dots,0)$.

\begin{prop} The group $H_m$ is a semidirect product of $L$ and $K$, where $L\unlhd G$ is unipotent and the subgroup $K\cong
{\rm\textbf{G}}_m$ is a maximal torus in $H_m$.
\end{prop}

We see that, regardless of the multiplicity of the exceptional
vertex, $T\cong \textbf{G}_m$ is a maximal torus in ${\rm Out}^K
(A_{\Gamma})$. From this we deduce the following theorem.

\begin{te}
Let $\Gamma$ be an arbitrary Brauer tree and let $A_{\Gamma}$ be a
basic Brauer tree algebra whose tree is $\Gamma$. Up to graded
Morita equivalence and rescaling there is unique grading on
$A_{\Gamma}$.
\end{te}
\noindent{\bf Proof.} By Proposition \ref{cokarak} homomorphisms
of algebraic groups from ${\textbf{G} }_m$ to ${ \textbf{G}}_m$
give us all gradings on $A_{\Gamma}$ up to graded Morita
equivalence. Since the only homomorphisms from ${ \textbf{G}}_m$
to ${ \textbf{G}}_m$ are given by maps $x\mapsto x^r$, for $x\in {
\textbf{G}}_m$ and $r\in \mathbb{Z}$, we have that there is unique
grading on $A_{\Gamma}$ up to rescaling (dividing each degree by
the same integer) and graded Morita equivalence (shifting each
summand of the tilting complex given by Green's walk).
$\blacksquare$

\section*{Acknowledgements}
I would like to thank Raphael Rouquier for suggesting the
problem of transfer of gradings via derived equivalence, and for
the constant support while I was writing this paper. Also, I would
like to thank Karin Erdmann for careful reading of the manuscript
and for many useful suggestions.

\end{document}